\documentclass[11pt]{article}
\usepackage{graphicx}
\usepackage{fancybox}
\usepackage{amsmath}
\usepackage{amssymb}
\usepackage{dsfont}
\usepackage[latin1]{inputenc}
\usepackage[T1]{fontenc}
\usepackage{textcomp}
\usepackage[a4paper]{geometry}
\usepackage{verbatim}
\usepackage{multirow}
\DeclareGraphicsExtensions{.jpg, .bmp}
\usepackage[pdfstartview=FitH]{hyperref} 
\usepackage[english]{babel}
\RequirePackage{natbib}

\def \M {\mathcal{S}^{d-1}}
\def \R {\mathbb{R}}
\def \L {\Lambda_{d-1}}
\def \P {\mathbb{P}}

\def \Rd {\mathbb{R}^{d}}

\def \a {\alpha}

\author{Antoine Dematteo, St\'{e}phan Cl\'{e}men\c{c}on, Nicolas Vayatis, Mathilde Mougeot}

\title{Sloshing in the LNG shipping industry: risk modelling through multivariate heavy-tail analysis}

\begin{document}
\maketitle
\begin{abstract}
In the liquefied natural gas (LNG) shipping industry, the phenomenon of \textit{sloshing} can lead to the occurrence of very high pressures in the tanks of the vessel. The issue of modelling or estimating the probability of the simultaneous occurrence of such extremal pressures is now crucial from the risk assessment point of view. In this paper, heavy-tail modelling, widely used as a conservative approach to risk assessment and corresponding to a worst-case risk analysis, is applied to the study of sloshing. Multivariate heavy-tailed distributions are considered, with Sloshing pressures investigated by means of small-scale replica tanks instrumented with $d\geq 1$ sensors. When attempting to fit such nonparametric statistical models, one naturally faces computational issues inherent in the phenomenon of dimensionality. The primary purpose of this article is to overcome this barrier by introducing a novel methodology. For $d$-dimensional heavy-tailed distributions, the structure of extremal dependence is entirely characterised by the angular measure, a positive measure on the intersection of a sphere with the positive orthant in $\mathbb{R}^d$. As $d$ increases, the mutual extremal dependence between variables becomes difficult to assess. Based on a spectral clustering approach, we show here how a low dimensional approximation to the angular measure may be found. The nonparametric method proposed for model sloshing has been successfully applied to pressure data. The parsimonious representation thus obtained proves to be very convenient for the simulation of multivariate heavy-tailed distributions, allowing for the implementation of Monte-Carlo simulation schemes in estimating the probability of failure. Besides confirming its performance on artificial data, the methodology has been implemented on a real data set specifically collected for risk assessment of sloshing in the LNG shipping industry.\\


\end{abstract}

\section{Industrial context}
In the liquefied natural gas (LNG) shipping industry, \textit{sloshing} refers to an hydraulic phenomenon which arises when the cargo is set in motion, \cite{PDS}. Following incidents experienced by the ships Larbi Ben M'Hidi and more recently by Catalunya Spirit, these being two LNG carriers faced with severe sloshing phenomena, rigorous risk assessments have become a strong requirement for designers, certification organisations (seaworthiness) and ship owners. In addition, sloshing has also been a topic of interest in other industries (for instance, see \cite{NASA} for a contribution in the field of aerospace engineering).
Gaztransport  $\&$  Technigaz (GTT) is a French company which designs the most widely used cargo containment system (CCS) for conveying LNG, namely the membrane containment system. The technology developed by GTT uses the hull structure of the vessel itself: the tanks are effectively part of the ship. The gas in the cargo is liquefied and kept at a very low temperature $(-163^{\circ}\mathrm{C})$ and atmospheric pressure, thanks to a thermal insulation system which prevents the LNG from evaporating. Although this technology is highly reliable, it can be susceptible to sloshing: waves of LNG apply very high pressures (over 20 bar) on the tank walls on impact and may possibly damage the CCS. Due to its high complexity, the sloshing phenomenon is modelled as a random process. The phenomenon is being studied by GTT experimentally on instrumented small-scale replica tanks ($1/40$ scale). The tanks are shaken by a jack system to reproduce the motion of the ship and induce the occurrence of sloshing, with the associated high pressures being recorded by the sensors. These experiments provide massive data sets which should hopefully, if adequately modelled, provide a better understanding of the spatial distribution of the pressure peaks and the dependence between them. As the tank is only partially instrumented, the structure of the dependence between extreme pressure values can only be observed locally in the tank where the sensors are installed. The next challenging step is to extrapolate the pressure field all around the tank based solely on the partial measurements provided by the sensors. This issue is not considered in the present article and will be the subject of further research (see \citet{Steinkohl11,Steinkohl12} for recent results on extreme value theory in the context of spatial models).

The phenomenon to be analysed here is described by a series of pressure measurements, and in particular by the existence of very large values corresponding to pressures created by heavy impacts, namely sloshing. Hence heavy-tail modelling is relevant in this context and is indeed considered as a conservative risk assessment method, insofar as it does not underestimate the importance of extreme values in general. Heavy-tailed distributions are also used for risk assessment in many other fields such as in finance \citep{RFM}, insurance \citep{mikosch} or for modelling natural hazards (refer to \cite{Tawn92} or \citet{Coles94}).

Modelling the one-dimensional marginal distribution of extreme observations is now common practice using the block maxima approach and the Generalized Extreme Value distribution (GEV), or the Peak Over Threshold approach (POT) and the Generalised Pareto Distribution (GPD) \citep{BGST,Reiss2007,Pickands_1975,Balkema.deHaan1974}. In contrast, the analysis of multivariate extreme data sets is much more challenging and this is the issue tackled in this paper. Since major damage occurs when the liquefied gas gives a heavy impact to a large area of the tanks, it is crucial to assess accurately the probability of simultaneous occurrences of very high pressures at several sensor locations. This paper considers the problem of estimating this key information. So far as the asymptotic study of multidimensional data sets is concerned, the vast majority of the results documented in the literature are mostly related to extreme-value parametric models \citep{Kluppelberg06,Boldi07}. Purely non-parametric approaches have also been considered, but their applications are generally restricted to the bivariate case \citep{Einmahl98,Einmahl2009}. Research into multivariate generalisations of the POT approach started only recently with the introduction of multivariate GPD, or even Generalized Pareto processes \citep{Buishand2008,Rootzen2006}. A few related simulation methods, limited to very specific models, are available.

In this paper, we develop a framework for accurately estimating the probability of failure of the containment system of LNG carrier tanks. A Monte-Carlo simulation scheme should ideally allow this probability to be approximated numerically. The target pressures are large, typically beyond the range of observed data. Assuming that sloshing data are derived from a multidimensional heavy-tail model, then when expressed in polar coordinates the radial part is asymptotically distributed as a generalized Pareto variable and independent of the angular component. The (asymptotic) distribution of the angular component is referred to as the \textit{angular measure} on the intersection of the unit sphere with the positive orthant of $Rd$. The extremal dependence between all $d$ sensors in the tank (or all the sensors in a specific area of the tank) is completely characterised by the angular measure. While simulation of the radius is straightforward using GPD distributions, simulating angles is challenging. When the tank is fitted with $d$ sensors, the angular measure can be decomposed into a mixture of up to $2^d$ - 1 sub-angular measures, with dimensions ranging from 0 to 
$d - 1$. Hence, any direct method for estimating the angular measure would suffer from the curse of dimensionality. phenomenon. An accurate understanding of the structure of the angular measure, that is of the asymptotic dependences between the sensors, is thus critical. Indeed, extremal pressures do not occur at the same time at all the different sensor locations and some sensors are likely to be asymptotically independent from some others. Hence, we seek a segmentation of the collection of sensors into $l$ groups such that: $(i)$ the measurements collected by the sensors in each subgroup are mutually independent in the extremes within this subgroup, $(ii)$ these measurements are mutually independent in the extremes from the other subgroups. Ideally, the cardinalities of the groups should be small with respect to $d - 1$ and $l$ small with respect to 2$^d-1$, so that estimation of the angular measure becomes tractable. For this purpose, we introduce here a novel methodology grouping the sensors into clusters satisfying assumptions $(i)$ and $(ii)$. This method is based on a spectral clustering algorithm \citep{Vonluxburg2007}, tuned to detect asymptotic dependences and independences. Ultimately, by conditioning upon the membership in each cluster, the asymptotic distribution of the data can be simulated and the corresponding risk of failure assessed.

The remainder of the paper is organized as follows. In section \ref{data}, we describe the data under study, explain how they have been collected and assess the relevance of heavy-tail modelling in the sloshing context. In section \ref{dependence}, we recall some basic concepts on multivariate regular variations and heavy-tail modelling extensively used throughout the paper. The method to perform the spectral clustering algorithm  tailored for multivariate extremes is presented in section \ref{spectralClustering}, and, based on the latter, estimation of the angular measure related to high dimensional observations is considered. In section \ref{caseStudy}, the technique promoted is next applied on real data in order to estimate the probability of simultaneous occurrence of high pressures in the tanks of LNG carriers and assess the risk induced by sloshing. In section \ref{discussion}, our main findings are discussed and possible lines of further research are sketched.

The remainder of the paper is organized as follows. In section \ref{data}, we describe the data under investigation, describe how they have been collected and assess the relevance of heavy-tail modelling in the context of sloshing. In section \ref{dependence}, we recall some of the basic concepts on multivariate regular variations and heavy-tail modelling that are extensively used throughout the paper. The methodology of the spectral clustering algorithm tailored for multivariate extremes is presented in section \ref{spectralClustering}, and, based on this method, estimation of the angular measure related to higher dimensional observations is considered. In section \ref{caseStudy}, the proposed technique is next applied to real data to estimate the probability of the simultaneous occurrence of high pressures in the tanks of LNG carriers and hence assess the risk induced by sloshing. In section \ref{discussion}, our main findings are discussed and possible lines of further research are outlined.

\section{Sloshing data and evidence of heavy-tail behaviour in sloshing events}
\label{data}
We start with a description of the sloshing data on which the subsequent statistical analysis relies, and then briefly review the basic concepts of heavy-tail analysis in extreme value theory, which have proved to be very relevant in the present context.

\subsection{Data set}
\begin{figure}[!htb]
\centering
\includegraphics[width=90mm,height=55mm]{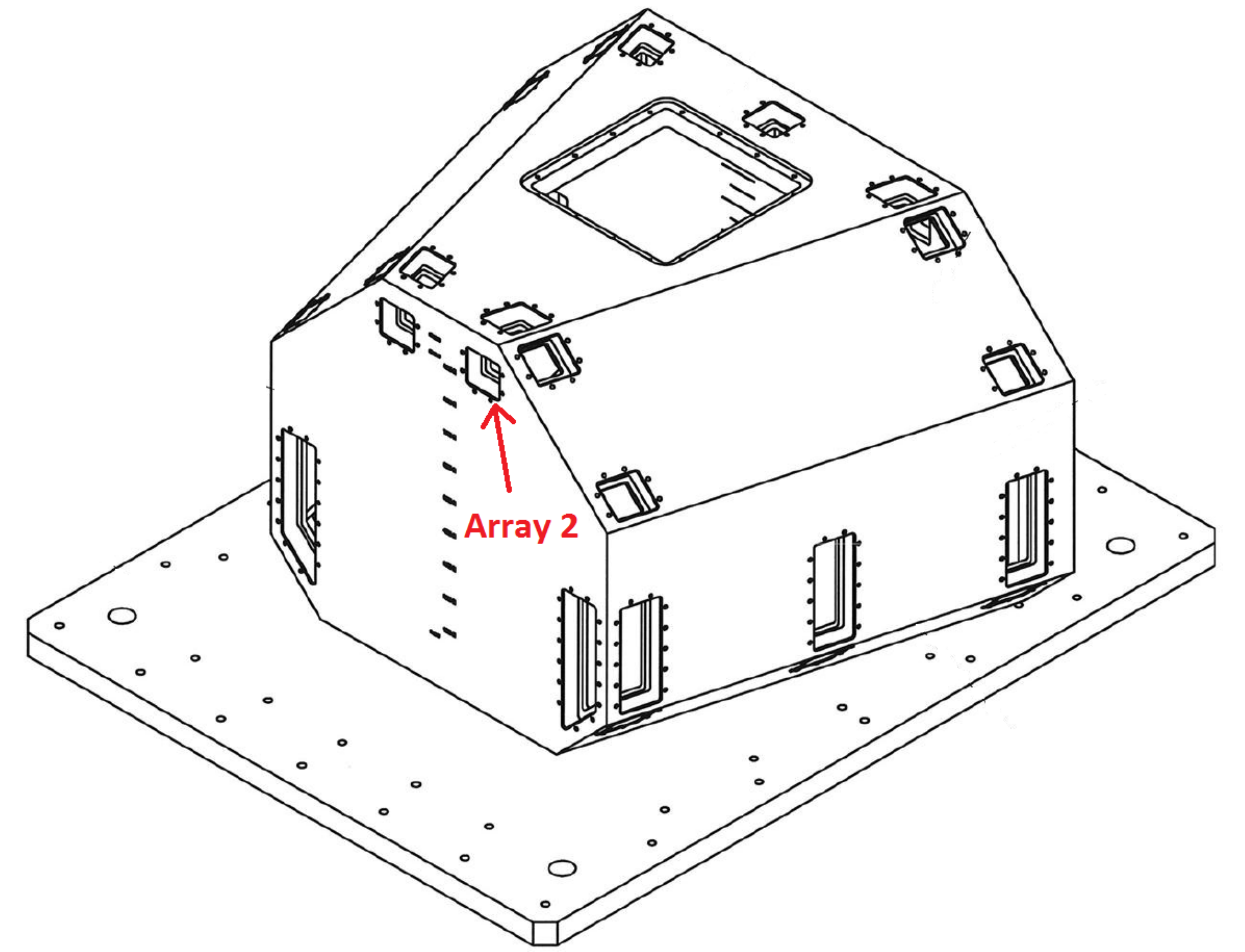}

\caption{. Diagram of a small-scale tank. The empty compartments are where the sensors are nested. We focus on the highly filled configuration in which the sensors measure the pressures recorded at the top of the tank.}
\label{small_tank}
\end{figure}

The data we consider here were provided by GTT and obtained during a test programme on small scale tanks (1/40 scale) as depicted in Fig. \ref{small_tank}. The small tank is filled with water (modelling the LNG) and SF6 gas (modelling the gaseous mixture lying above the LNG in the tank). Here, the density ratio between SF6 gas and water is the same as that between LNG and the mixture \citep{Brosset2009}. The tank replica is shaken by a jack system to reproduce the ship motions. The tank is instrumented with a collection of sensors grouped into arrays. As soon as a sensor records a pressure above a threshold, the pressures measured simultaneously by all the other sensors of the array are recorded also at a sampling frequency of 20kHz until the pressure signal falls below the threshold for each sensor. The signal recorded by a sensor after a pressure peak exhibits a typical sinusoidal shape and decreases slowly. In this study, for each high pressure event and for each sensor, risk assessment is based on the pressure peak only.

\begin{figure}
\centering
\makebox{\includegraphics[width=40mm,height=40mm]{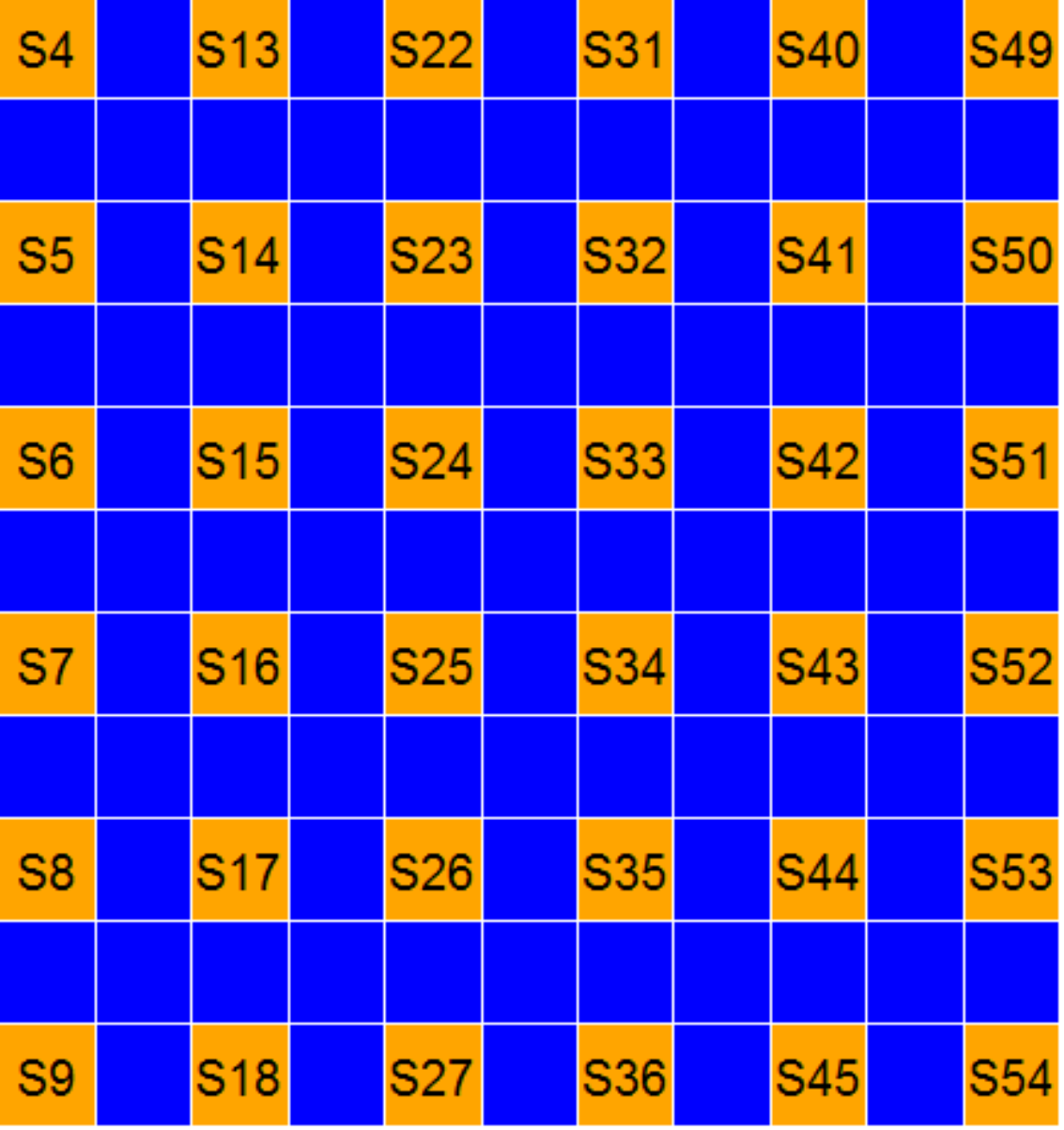}}
\caption{\label{module}Array number 2, a $6\times 6$ sensors array.}
\end{figure}

The data set provided by GTT corresponds to a \textit{high filling configuration} where the tanks are nearly full of liquefied gas. A diagram of the small-scale tank is shown in Fig. \ref{small_tank}. We focus on array number 2 (see Fig. \ref{module}), with $d = 36$ sensors on this array. The total number of raw observations per sensor is $n = 145,326$, which corresponds to 6 months of sailing at full scale. Fig. \ref{nb_evt_sens} in appendix \ref{statdesc} shows a map of the number of impacts detected at different locations of the array.

\begin{figure}[!htb]
\centering
\includegraphics[width=145mm,height=65mm]{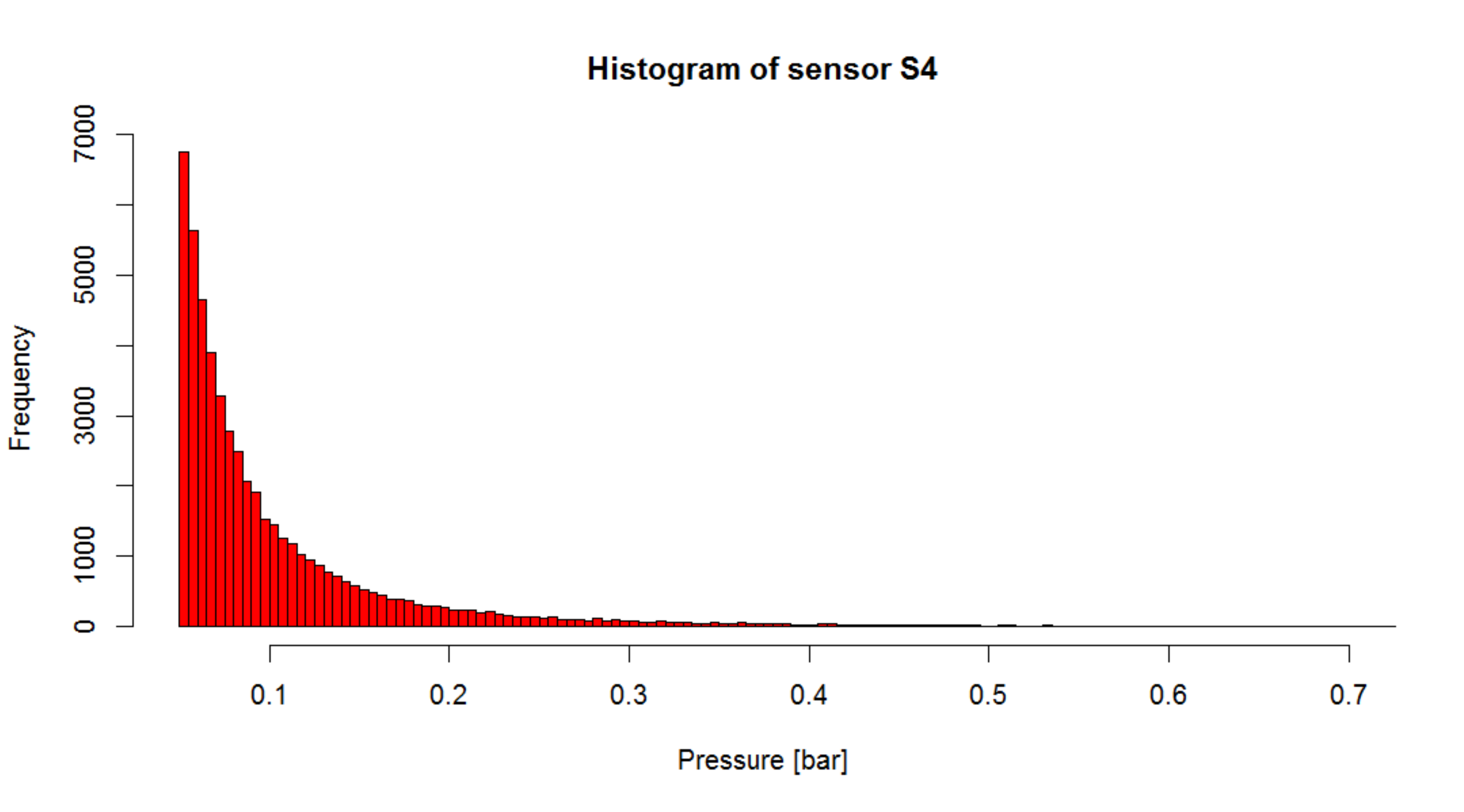}

\caption{Histogram of pressure measurements for sensor $S4$. Only pressures smaller than 0.7 bar are shown.}
\label{histogramS4}
\end{figure}

In the high filling configuration, a pressure measurement is considered as a sloshing impact if it is above 0.05 bar. For example, sensor $S4$ records $52769$ such impacts. Fig. \ref{histogramS4} is a histogram of the pressure values for this sensor. This histogram, together with Table \ref{quantiles}, shows that, even a long way from the mean, many high pressure events can be observed and this gives us a first insight into the clear relevance of heavy-tail modelling in this context. More detailed statistics are provided in Table \ref{quantilesallsens} in appendix \ref{statdesc}

\begin{table}
	\caption{\label{quantiles}Extreme quantiles of sensor $S4$. The maximum observed is 1.74 }
\fbox{%
\begin{tabular}{cccccc}
\textbf{\textit{order}}       & \textbf{0.9}   & \textbf{0.99} & \textbf{0.999} & \textbf{0.9996} & \textbf{0.9999} \\
\textbf{\textit{value} [bar]} & 0.19  & 0.48 & 0.87  & 1.00   & 1.20   \\
\end{tabular}}
\end{table}

\subsection{Heavy-tail analysis}
\label{HTanalysis}

By definition, heavy-tail phenomena are those which are governed by very large values, occurring with a non-negligible probability and with significant impact on the system under study. When the phenomenon of interest is described by the distribution of a univariate random variable, the theory of regularly varying functions provides the appropriate mathematical framework for heavy-tail analysis. For the sake of clarity, and in order to introduce some notation to be used in what follows, we recall some related theoretical background. Refer to \citet{resnick2007heavy}, \citet{Hult05} and \citet{Hult06} for an account of the theory of regularly varying functions and its application to heavy-tail analysis.

Let $\a>0$. We denote by 
$$
\mathcal{RV}_{-\a}=\{U:\R^+\rightarrow\R^+\text{ Borel measurable } |\lim_{t\rightarrow\infty}\frac{U(tx)}{U(t)}=x^{-\a},\: x>0\}$$

the space of regularly varying functions with index $\a$. Let $X$ be a random variable with cumulative distribution function (\textit{cdf}) $F$ and survival function $\overline{F}=1-F$. The random variable $X$ is said to have a heavy (right) tail of index $\a$ when $\overline{F}\in\mathcal{RV}_{-\a}$. The cdf $F$ of any heavy-tailed random variable with tail-index $\a$ can be written  as $F(x)=1-L(x)\,x^{-\a}$, where $L$ is a slowly varying function, \textit{i.e.} $L\in\mathcal{RV}_{0}$. In addition, the heavy-tail property can be classically formulated in terms of vague convergence to a homogeneous positive measure. Indeed, the random variable $X$ belongs to $\mathcal{RV}_{-\a}$ if and only if:

$$
n\mathbb{P}\left(X/F^{-1}(1-1/n)\in \cdot  \right) \overset{v}{\rightarrow} \mu_{\alpha}(\cdot) \text{ in } M_+(0,\, \infty],
$$
where $F^{-1}(u)=\inf\{t:\; F(t)\geq u\}$ denotes the generalised inverse of $F$, $\mu_{\alpha}(x,\, \infty]=x^{-\alpha}$, $M_+(0,\, \infty]$ the set of non-negative Radon measures on $(0,\, \infty]$ and $\overset{v}{\rightarrow}$ vague convergence.

\begin{figure}
\centering
\makebox{\includegraphics[width=145mm,height=90mm]{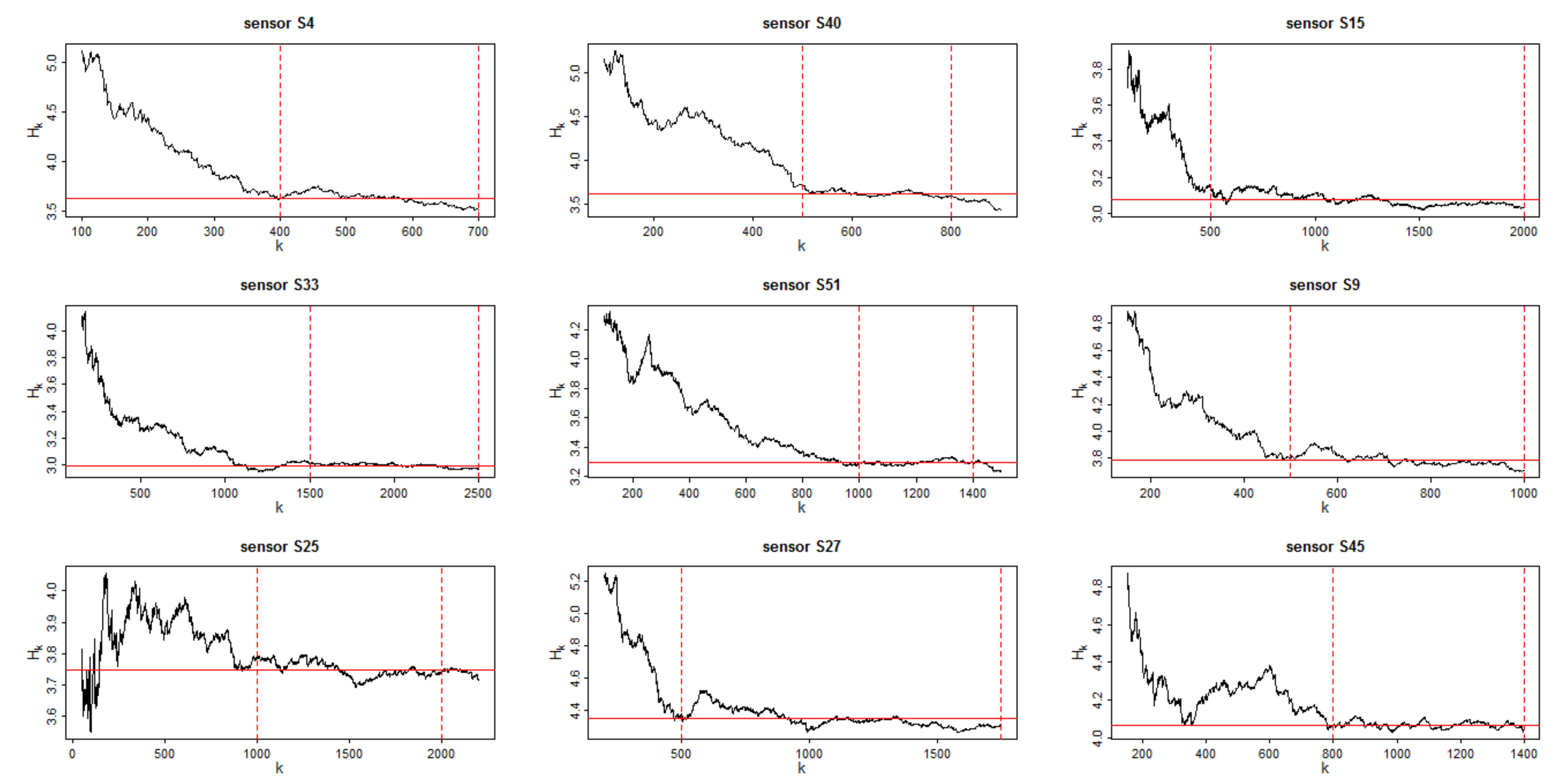}}
\caption{\label{Hillplot3}Horror Hill plots for 9 sensors in the module. The dotted vertical lines show the regions where the plots appear nearly constant. The solid horizontal line gives the estimated value of  $\alpha$.}
\end{figure}

The tail index $\a$ can be estimated by the popular Hill estimator, see \cite{Hill1975}:

$$
\widehat{\alpha}_{k,n}=\left(\frac{1}{k}\sum_{i=1}^{k}\text{log}\left(\frac{X(i)}{X_{(k)}}\right)\right)^{-1},
$$
where we denote by $X(1)>\dots>X(n)$ the order statistics of $X_1,\dots,X_n$.
This estimator is consistent and asymptotically normal under certain assumptions, \textit{i.e.}: $\sqrt{k}\left(\widehat{\alpha}_{k,n} - \alpha\right)$ converges in distribution to a centred Gaussian random variable with variance $\alpha$ as $k\rightarrow \infty$ such that $k=o(n)$. However, its behaviour  can be very erratic in $k$. In practice, to handle the possible variability in $k$, we plot the graph of the mapping $k\mapsto \widehat{\alpha}_{k,n}$ and seek a region where the resulting \textit{Horror Hill Plot} is nearly constant, see  for instance \citet[Chap.9]{resnick2007heavy}.  Fig. \ref{Hillplot3}shows the related Hill Horror plots. Table \ref{HillCI} gives estimates for $\a$ for all sensors in the array, and shows no evidence of variation of the tail index along the tank.

\begin{table}
\caption{\label{HillCI} Hill estimate of $\alpha$ for the sensors of the array and 90\% Gaussian confidence interval. The estimates of the table correspond to the locations of the sensors indicated in Table \ref{module}}
\fbox{%
\begin{tabular}{cccccc}
$\mathbf{3.63}$ &$\mathbf{3.50}$ &$\mathbf{3.16}$ &$\mathbf{3.51}$ &$\mathbf{3.61}$ &$\mathbf{3.79}$ \\
(3.27-3.99) & (3.24-3.75) & (2.83-3.49) & (2.95-4.07) & (3.29-3.92) & (3.44-4.14) \\
 $\mathbf{3.60}$ &$\mathbf{3.55}$ &$\mathbf{3.08}$ &$\mathbf{4.15}$ &$\mathbf{4.35}$ &$\mathbf{4.09}$ \\
(3.25-3.96) & (3.03-4.08) & (2.73-3.42) & (3.83-4.46) & (4.02-4.67) & (3.77-4.42) \\
 $\mathbf{4.12}$ &$\mathbf{3.62}$ &$\mathbf{3.51}$ &$\mathbf{3.75}$ &$\mathbf{4.32}$ &$\mathbf{4.35}$ \\
(3.86-4.37) & (3.30-3.93) & (2.95-4.07) & (3.56-3.94) & (3.87-4.77) & (4.08-4.62) \\
 $\mathbf{4.30}$ &$\mathbf{3.26}$ &$\mathbf{3.00}$ &$\mathbf{3.60}$ &$\mathbf{3.95}$ &$\mathbf{4.12}$ \\
(4.01-4.59) & (2.81-3.71) & (2.81-3.18) & (3.38-3.82) & (3.67-4.24) & (3.78-4.47) \\
 $\mathbf{3.62}$ &$\mathbf{3.25}$ &$\mathbf{3.63}$ &$\mathbf{3.85}$ &$\mathbf{4.44}$ &$\mathbf{4.07}$ \\
(3.38-3.86) & (3.05-3.45) & (3.33-3.94) & (3.47-4.22) & (4.00-4.87) & (3.67-4.47) \\
 $\mathbf{3.65}$ &$\mathbf{3.35}$ &$\mathbf{3.29}$ &$\mathbf{3.90}$ &$\mathbf{4.28}$ &$\mathbf{4.05}$ \\
(3.33-3.97) & (3.13-3.57) & (3.01-3.58) & (3.62-4.18) & (3.98-4.57) & (3.72-4.37) \\
 
\end{tabular}}
\end{table}

\section{Dependency in the extremes : multivariate heavy-tail analysis and angular measure}
\label{dependence}

Based on a sample of $i.i.d.$ observations $X_1, \dots, X_n$, the aim is now to investigate the dependence structure among the large pressures simultaneously measured by different sensors and to implement adequate statistical methods in order to reliably estimate the probability that several sensors simultaneously record extremal pressures (possibly outside the range of the data sample). It should be emphasised that the focus is on observations far from the mean behaviour. Note that simple moment-based quantities such as covariance matrices are clearly inadequate for describing dependences on extremal values. Indeed they do not distinguish between dependence among large or small values, which might rely on very different mechanisms. For multivariate heavy-tailed data, as recalled below, the dependence structure in regard to extremal observations is fully described by the notion of \textit{angular measure}.

\subsection{Notations}
\label{Notations}

\paragraph{General Notation:\\}
Here and throughout, we consider a collection $X = (X^{(1)},\; \ldots ,\;  X^{(d)})$  of pressures , drawn from a probability distribution $F(dx)$, measured by a group of $d \geq 2$ sensors: $X^{(i)}$ is the pressure measured by sensor number $i$ and $F_i(dx)$ denotes its marginal probability distribution.  The cumulative distribution function of the random variable  $X$ is given by $F(\mathbf{t}) = \P\{X^{(1)}<t^{(1)},\;\ldots,\; X^{(d)}<t^{(d)}\}$ for all $\mathbf{t} = (t^{(1)},\; \ldots,\; t^{(d)})\in \mathbb{R}^d_{+}$. Finally we denote by $\mathbf{u} = (u, \dots , u)$  the $d$-dimensional vector whose coordinates are all equal to $u\in \bar{\mathbb{R}}$ and by $u \cdot\mathbf{t}$ the vector $(u\cdot t^{(1)}, \dots , u \cdot t^{(d)})$. In addition, all operations in what follows are taken to be component-wise and for $t \in \R^{+}$, $X>t$ means that all the components of the vector $X$ are greater than $t$. 

\paragraph{Standardisation:\\}
We denote by $Z=\left(Z^{(1)},\dots, Z^{(d)}\right)$ the random variable whose components are given by
 \begin{equation}\label{eq:stand}
 Z^{(i)}=1/(1-F_i(X^{(i)})),\: i = 1,\; \ldots,\;  d,
 \end{equation} so that each margin of the vector $Z$ is standard Pareto distributed, \textit{i.e.} $\P(Z^{(i)}>x)=1/x,\: i=1,\; \ldots,\; d$. In practice, as the $F_i$'s are unknown, they may be replaced by their empirical counterparts in \eqref{eq:stand}. This technique, used in the subsequent analysis, is referred to as the \textit{ranks method} (see \citet[][subsection 9.2.3]{resnick2007heavy} for further details).

\paragraph{Set notations and specific sets:\\}
The indicator function of any event $\mathcal{E}$ is denoted by $\mathds{1}(\mathcal{E})$. The Dirac measure associated with any set $A$ is denoted by $\delta_{A}$ and its complementary subset by $A^{c}$. The punctured positive orthant is denoted by $\mathcal{O}=\mathbb{R}^{d}_+\backslash\{0\}$.  For a given norm $\|.\|$ on $\mathcal{O}$,  the set $\L$ is the intersection of the unit sphere (with respect to the chosen norm) $\M := \{x \in \Rd ,\|x\| = 1\}$ with $\mathcal{O}$.

The norms defined by $\vert\vert x\vert\vert_p=\left(\sum_{j=1}^d\vert x_j\vert^{p}\right)^{1/p}$ and  $\vert\vert x\vert\vert_{\infty}=\max_{i=1\dots d} \vert x_i\vert$ for all $x=(x_1,\;\ldots,\; x_d)\in\mathbb{R}^d$ are referred to as the $\mathcal{L}_{p}$-norm and $\mathcal{L}_{\infty}$-norm.

The set of all partitions of $\{1,\dots,d\}$ is denoted by $\mathcal{P}_d$. For an element $p=(i_1,\dots,i_m)\in \mathcal{P}_d$, we denote by $\bar{p}=\{1,\dots,d\}\backslash p$ and by $X^{(p)}=\left(X^{(i_1)},\dots,X^{(i_m)}\right)$. The number of elements in $p$ is denoted by card($p$).

\subsection{Standard case: identical tail index for all sensors}
\label{standard}

In the standard case, all marginal distributions are tail equivalent, meaning that they have the same index $\alpha = \alpha_1 = \ldots =\alpha_d.$ In this case, the probability distribution $F(dx)$ is said to be regularly varying with index $\alpha$ when there exists a Radon measure $\nu(dx)$ on  $\mathcal{O}$ such that

\begin{equation}
\label{regvarmu}
\lim_{\lambda \rightarrow \infty} \frac{1-F(\lambda \mathbf{t})}{1-F(\lambda)}=\nu\left(\left[0,\mathbf{t}\right]^{c}\right),
\end{equation}
the measure $\nu$ having the homogeneity property : $\nu\left(\left[0,u.\mathbf{t}\right]^{c}\right)=u^{-\alpha}\times\nu\left(\left[0,\mathbf{t}\right]^{c}\right)$.

Multivariate heavy-tailed distributions are conveniently described using polar coordinates. Consider two norms s $\|.\|_{(1)}$ and $\|.\|_{(2)}$ on $\mathbb{R}^{d}$ and define
$T:x\in \mathcal{O} \mapsto ( \|x\|_{(1)},x/\|x\|_{(2)})\in \mathbb{R}^{\star}_{+} \times \L$.
 For notational simplicity, we set $(r, a) = T(x)$ as well as $(R,A) = T (X)$ when considering random variables. Condition \eqref{regvarmu} can be then formulated as follows: there exists a constant $c \in \mathbb{R}_{+}$ and a probability measure $S(da)$ on $\L$ such that,
\begin{equation}
\label{limpolar}
u\P\Bigg\{\left(\frac{R}{b(u)},A\right)\in\left[0,r\right]^{c}\times \Theta\Bigg\}\xrightarrow[u\rightarrow\infty]{} c \cdot r^{-\alpha} \times S(\Theta)~:=\nu\circ T^{-1}\left( \left[0,r\right]^{c}\times \Theta\right),
\end{equation}
for any Borel set $\Theta \subset \L$, any $r > 0$.  The function $b(u) = F_R^{-1}(1-\frac{1}{u})$ is the $(1-\frac{1}{u})$-quantile of the distribution of $R$.   The limiting measure $\nu$ is referred to as the \textit{exponent measure}. The measure $S$ is known as the \textit{angular measure} and provides a complete description of the tail dependence structure. When concentrated around the intersection of the  line $\left\{x\in \R^{d} : x_1=\dots=x_d\right\}$ and $\L$ (the point of coordinates $(0.5,0.5)$ in the bivariate case when considering the $\mathcal{L}^{1}$-norm), a tendency toward complete extremal dependence can be observed. In contrast, if the angular distribution is concentrated at the intersection of $e_j$ with $\L$, $1\leq j\leq d$, where $e_j$ is the unit vector with coordinates $0$ everywhere except along the $j$'th axis, then there is a tendency towards complete independence.

A natural estimator of the angular measure is defined as follows. Set a large threshold $t$ and apply the polar operator to the rank transformed data $Z$ to obtain $\big((R_i,A_i),\:i=1\dots n\big)$. The estimate $\widehat{S}$ of $S$ is:

\begin{equation}
\label{angularEst}
\widehat{S}(\Theta)=\sum_{i=1}^{n}\mathds{1}\left(A_i\in \Theta, R_i>t\right)
\end{equation}

The estimated angular measure $\widehat{S}(\Theta)$ can be normalised by $\widehat{S}(\L)$ to become the probability distribution $\widehat{S}(\Theta)/\widehat{S}(\L)$. For simplicity, throughout the paper we shall continue to denote by $\widehat{S}(\Theta)$ the angular probability measure. When attempting to estimate directly the density of the (supposedly absolutely continuous) angular probability by means of kernel smoothing techniques for instance, we may face major computational difficulties inherent in\textit{ the curse of dimensionality}, even for moderate values of dimension $d$. As shown in the previous section, heavy-tail modelling is quite appropriate in the context of sloshing data. However, it needs to be combined with an adequate dimension reduction technique before carrying out any statistical procedure. 

\subsection{Decomposition of the angular measure}
\label{mixtureModel}
\begin{figure}[!htb]
\centering
\includegraphics[width=145mm,height=90mm]{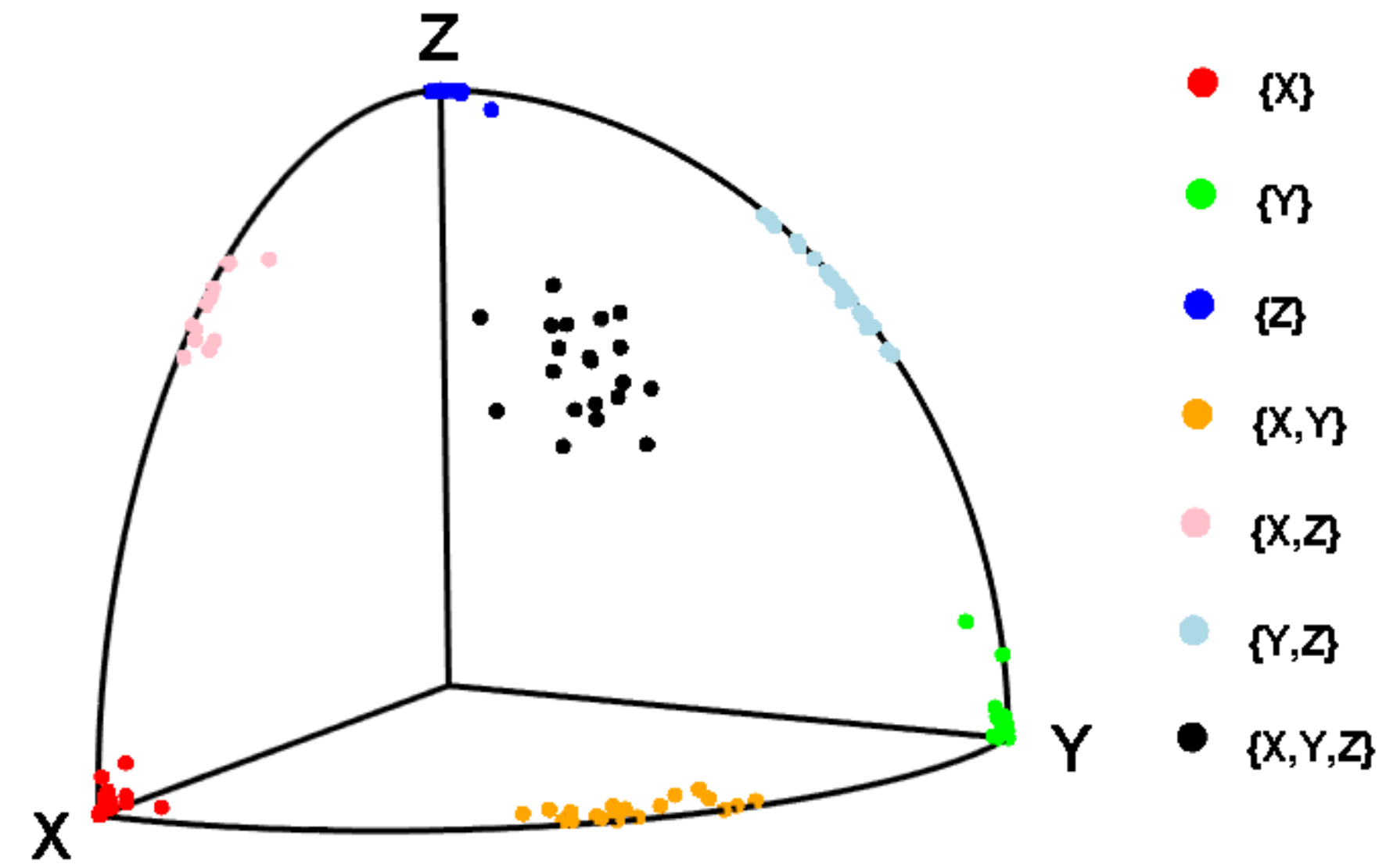}

\caption{Extreme data points projected on $\L$. The data have been simulated so that all the faces are support faces. The asymptotic dependence for each point is indicated on the right. Clusters can be easily identified}
\label{3dexemple}
\end{figure}

In the subsequent analysis, we denote the angular probability measure by $S:=S/S(\L)$ . The extreme dependence structure between $d$ variables $(X_1,\dots,X_d)$ is entirely characterised by the angular probability $S$ and more specifically by the geometry of its support, denoted by $\textit{supp}(S)$ and included in the set $\L$. This set is the reunion of $2^d-1$ open faces of dimensions ranging from $0$ to $d-1$. Denote the set of all these faces by $\mathcal{F}_d$ . There is a one-to-one correspondence between $\mathcal{P}_d$ and $\mathcal{F}_d$ and we have $\textit{supp}(S)\subset \mathcal{F}_d$. More precisely, for any element $p_m=\{i_1,\dots, i_m\}\in \mathcal{P}_d$, with $1\leq m \leq d$, if the variables $X_{i_1},\; \ldots,\; X_{i_m}$ exhibit asymptotic dependence, the support of their (sub-)angular probability $S_{p_m}$ is non empty and has dimension $m-1$.By contrast, in the case of asymptotic independence, the support of the angular measure is empty. These considerations suggest the following \textit{mixture model} for the angular probability distribution:

\begin{equation}
\label{mixtureangprob}
S=\sum_{p\in \mathcal{P}_d}\pi_p S_{p},
\end{equation}
where $\sum_{p\in \mathcal{P}_d}\pi_p = 1$ and for any $p=\{i_1,\dots,i_m\}$ with $1\leq m\leq d$, $\pi_p=S(supp(S_p))$, \textit{i.e.} it is the proportion of observations for which the variables $X_{i_1},\dots,X_{i_m}$ are jointly extreme.
The angular components of the largest (polar transformed) observations form clusters of points on $\L$, each cluster being contained in a face of $\mathcal{F}_d$ (Fig. \ref{3dexemple} provides a simulated example in dimension 3). In order to characterise the dependence structure of $(X_1,\;\ldots,\; X_d)$, we need to identify the sub-angular measures $(S_p)_{p\in \mathcal{P}_d}$ with non empty supports which boils down to identifying the clusters or the associated support faces of $\mathcal{F}_d$. The methodology for achieving this aim is introduced in the next section and is inspired by spectral clustering techniques.

\section{Spectral Clustering : recovering the faces}
\label{spectralClustering}
The purpose of the methodology introduced in this section is to provide a sound estimate of the angular probability in the presence of large-dimensional data sets. In subsection \ref{mixtureModel}, we introduced a mixture model explicitly stating that the angular probability is a weighted sum of up to possibly $2^{d}-1$ sub-angular probabilities  with dimensions ranging from 0 to $d-1$. Based on a theoretical framework introduced in subsection \ref{spectClust}, the aim of subsection \ref{asymptdepest} is to identify the sub-angular probabilities that are not identically zero. Assuming there are $l$ such probabilities with dimension $d(1),\dots,d(l)$, the dimensionality will be efficiently handled by the decomposition of Eq. \eqref{mixtureangprob}  if the two following conditions hold:
\begin{itemize}
\item[(i)] The number $l$ is small with respect to $2^{d}-1$ so that there are not too many terms in the sum of Eq. \eqref{mixtureangprob}.
\item[(ii)] The maximal dimension $\max_{i=1\dots l}d(i)$ is small with respect to $d$. 
\end{itemize}

If these two conditions are satisfied, the estimation of the angular probability of $(X^{(1)},\dots,X^{(d)})$ will be tractable.
 
\subsection{Theoretical background to spectral clustering.}
\label{spectClust}

\textit{Spectral clustering} is a segmentation technique quite adapted to data lying on a Riemannian manifold since the metric used to describe the distance between data points can be chosen in a very flexible manner, see \citet{Vonluxburg2007}. In particular, a Riemannian metric on $\L$ can be considered for this purpose. In addition, a significant advantage of spectral clustering as compared with to certain alternative clustering techniques is that it does not require the number of clusters describing the data to be set in advance, i.e. in our case, the number of support faces. For clarity, we start off with recalling briefly the rationale behind the spectral clustering approach

Given a data set $(x_1, \dots, x_n)$ and coefficients $w_{i,j} \geq 0$ measuring the similarity between all pairs of observations $(x_i,x_j)$, we can construct a similarity graph $G=(V,E)$. Each vertex $v_i$ represents a data point $x_i$. Two vertices are connected if the similarity $w_{i,j}$ between the corresponding data points $x_i$ and $x_j$ is strictly positive and the edge is weighted by $w_{i,j}$. The clustering algorithm aims to find a partition of the graph such that the similarities between vertices of a same cluster are greater than those between vertices lying in different groups.
A typical choice for quantifying similarity is the Gaussian function $w_{i,j}=\exp(-\rho_{i,j}^2/2\sigma^2)$, where the parameter $\sigma$ controls the width of the neighbourhoods and $\rho_{i,j}$ is the Riemannian distance between $x_i$ and $x_j$.
Some additional notations and definitions are required in order to describe the spectral clustering mechanism. The {\it weighted adjacency matrix} of the graph is $W = (w_{i,j})_{1\leq i,j\leq n}$. As the graph $G$ is undirected, we require $w_{i,j} = w_{j,i}$. The {\it degree} of a vertex $v_i \in V$ is defined as $d_i = \sum_{j=1}^{n} w_{i,j}$ and the {\it degree matrix} $D$ as the diagonal matrix with the degrees $d_1,\; \ldots,\; d_n$ as diagonal coefficients. 

\begin{figure}[!htb]
\centering
\includegraphics[width=130mm,height=100mm]{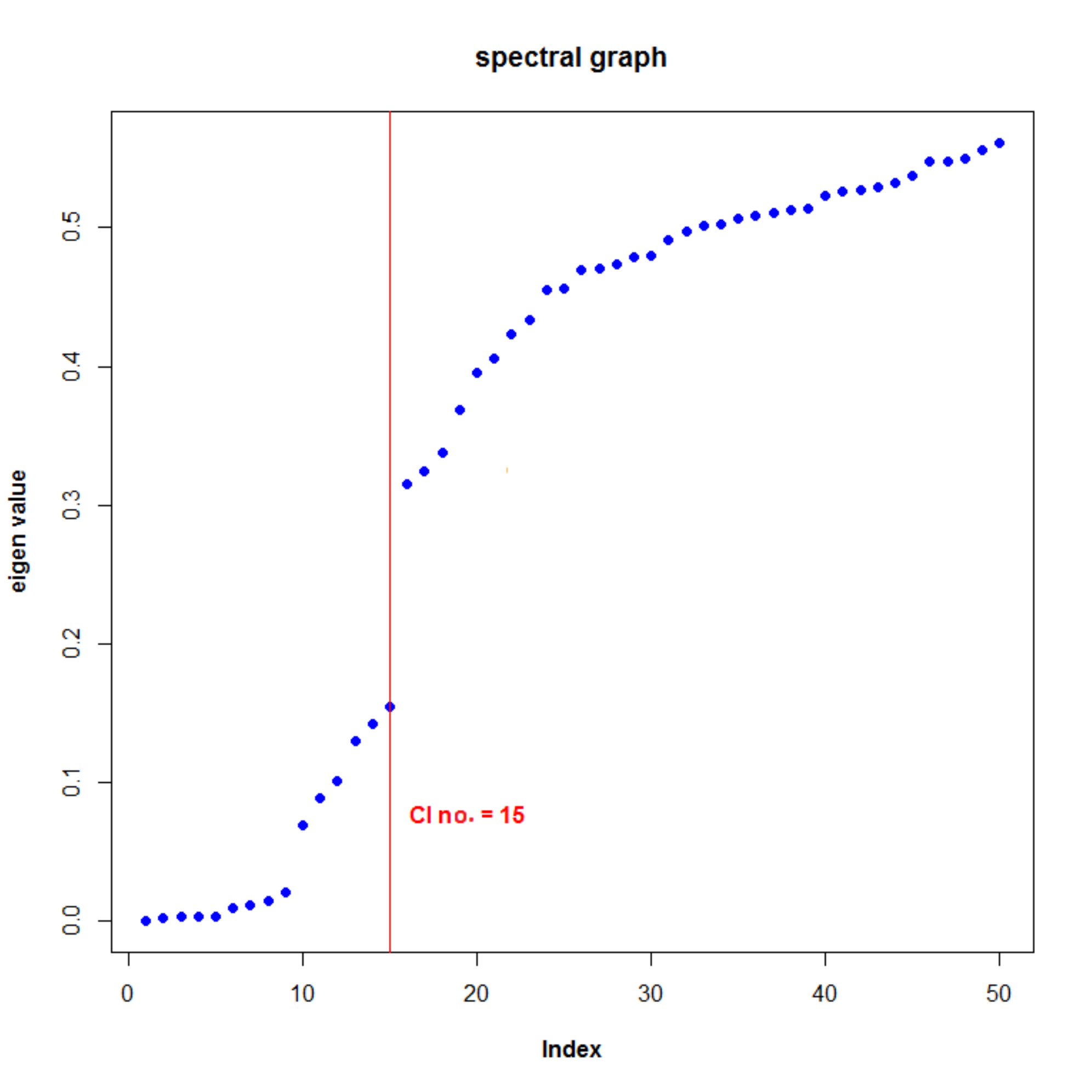}

\caption{Graph of the first 40 eigenvalues. For this simulation, a gap was detected between the 15th and the 16th eigenvalues, indicating 15 sets of asymptotically dependent variables.}
\label{SpectralGraphSimu}
\end{figure}

Armed with these notations, the {\it graph Laplacian} is $L= D - W$  and the {\it normalised graph Laplacian} is defined by  $L_{sym} =  D^{-1/2}LD^{-1/2}$. The matrix $L_{sym}$ exhibits some very interesting properties: the multiplicity $k$ of the eigenvalue $0$ of $L_{sym}$ is equal to the number of connected components $A_1,\; \ldots,\; A_k$ in the graph and the eigenspace corresponding to the eigenvalue $0$ is spanned by the related indicator vectors $\mathds{1}_{A_1} , \dots , \mathds{1}_{A_k}$ (in practice, the eigenvalues of $L_{sym}$ are not strictly zero and one needs to detect a gap. See Fig. \ref{SpectralGraphSimu} and \ref{SpectralGraph} for an illustration). Based on these results,  \citep{Ng2002} proposed the clustering algorithm presented in appendix \ref{SpectralClustering}. It involves the popular $k$-means vector quantization method (see \cite{Hartigan75}). We point out that the clustering produced by the $k$-means algorithm corresponds to a local optimum and depends strongly on the initialisation parameters. In practice the algorithm must therefore be run several times.

\medskip

\subsection{Application to asymptotic dependence estimation}
\label{asymptdepest}
In this section, we derive an algorithm for finding the groups of asymptotically dependent variables. We consider the standardised observations $Z_1,\dots,Z_n$ and apply the polar transform $(R_i,A_i)=T(Z_{i}),\: i=1\dots,n$. We consider the extreme data set $\Theta(t):=\left\{A_i\vert R_i>t, i=1\dots n\right\}$, where $t$ is a large threshold. 

The spectral clustering algorithm is used to infer the optimal number $l$ of clusters in the data set $\Theta(t)$ as well as the clusters $C_1,\dots,C_l$ themselves. The support face of each cluster $C_i,\: i = 1\dots l,$ is in a one-to-one correspondence with a group $E_i\in\mathcal{P}_d$ of asymptotically dependent variables; owing to some potential pitfalls, this needs to be estimated with care. The caveats associated with this estimation of $E_i$ are better understood via the concepts of coefficient of tail dependence $\eta$ (\citet{Ledford96}) or also by hidden regular variations (see subsection 9.4 in \citet{resnick2007heavy}). 

In practice, statistical methods may experience difficulties in distinguishing between asymptotic independence and exact independence, and also between asymptotic dependence and independence. For instance, if $\eta\rightarrow 1/2^{-}$, the variables are asymptotically independent but even for very large values they are likely to co-occur.

Based on these observations, we propose a heuristic technique for estimating $E_i$. Formally, with each cluster $C_i$ of size $c_i$, associate a threshold $e_i:=e_i(c_i)$ and define
$$E_i := \left\{j=1,\dots,d \:\left\vert\:\sum_{l \in C_i}\mathds{1}\left(Z^{(j)}_l>t\right)\geq e_i \right.\right\}.$$

The following \textit{extremal spectral clustering} algorithm is derived from the above considerations.

\begin{center}
\fbox{
\begin{minipage}{0.8\textwidth}
\medskip

\centering
\textbf{\textit{Extremal Spectral Clustering}}\\
\begin{itemize}
\item[]\textbf{Input}: i.i.d sample of size $n$ of $Z=\left(Z^{(1)},\dots,Z^{(d)}\right)$, standard Pareto distributed. 
\item[]\textbf{Parameters}: Threshold $t$. Number $n_r$ of repetition of the $k$-means algorithm. Minimal number $m_r$ of acceptance of a cluster. 

\item Apply the polar transform $(R_i,A_i)=T(Z_i),\: i=1,\dots,n$
\item Form the set $\Theta(t):=\left\{A_i\vert R_i>t, i=1\dots n\right\}$. Assume $card({\Theta(t)})=K$.
\item Compute $D\in\R^{K\times K}$ where $D_{i,j}$ is the Riemannian distance on $\L$ between $A_i$ and $A_j$.
\item Repeat the Spectral clustering algorithm $n_r$ times,  with similarity matrix $D$ as input. Select the clusters appearing at least $m_r$ times. Denote them by $C_1,\dots, C_l$, their size by $c_1,\dots,c_l$ and the thresholds by $e_1,\dots,e_l$
\item For any $i=1\dots l$ derive the set $E_i$ from $C_i$. Some $E_i$ might be empty and others might appear several times. Denote by $E_1,\dots,E_{l_0}$ the unique non empty sets.
\item[]\textbf{Output}: $E_1,\; \ldots ,\; E_{l_0}$.
\begin{flushleft}
\end{flushleft}
\end{itemize}

\end{minipage}
}
\end{center}

In the remaining of this paper, unless explicitly stated, $t$ is the threshold used to distinguish between extreme and non extreme observations and we set $\|.\|_{(1)}=\|.\|_{\infty}$ and $\|.\|_{(2)}=\|.\|_{2}$. In addition, $l$ will always stand for the number of support faces of the angular measure of $\mathbf{Z}$ and $E_i,\: i =1\dots l$ are the associated sets indexing the asymptotically dependent variables. \\

Once the groups $E_i,\: i=1\dots l$ have been estimated, estimation of each sub-angular probability (density respectively) is straightforward using the empirical estimate of Eq. \eqref{angularEst} (kernel estimators respectively) so that the only issue is the estimation of the coefficients $\pi_p,\:p\in\mathcal{P}_d$. We define the sets $\mathfrak{P}_d=\left\{E_i, \: i = 1\dots l\right\}$, which is a subspace of $\mathcal{P}_d$, and $ \mathcal{I}_t=\left\{ i ,\: \|Z_i\|_{\infty} > t \right\}$. We set $N_t=\text{card}(\mathcal{I}_t)$. The estimator for $\pi_p$ is defined as follows:

\[\pi_{p}=\left\{
   \begin{array}{lr}
		0 & \text{if } p \notin \mathfrak{P}_d\\
		\frac{1}{N_t} \sum_{i\in\mathcal{I}_t}\mathds{1}\left(Z_i^{(p)}>t , Z_i^{(\bar{p})}<t\right) & \text{otherwise}\\
	\end{array}\right.
.
\]

\subsection{Estimating the probability of joint exceedance}
\label{jointExcProb}
In what follows,  it is assumed that on each of the $l$ support faces, the angular measure has a density with respect to the Lebesgue measure on the associated support face. This density is referred to as angular density.

\begin{figure}[!htb]
\centering
\includegraphics[width=100mm,height=100mm]{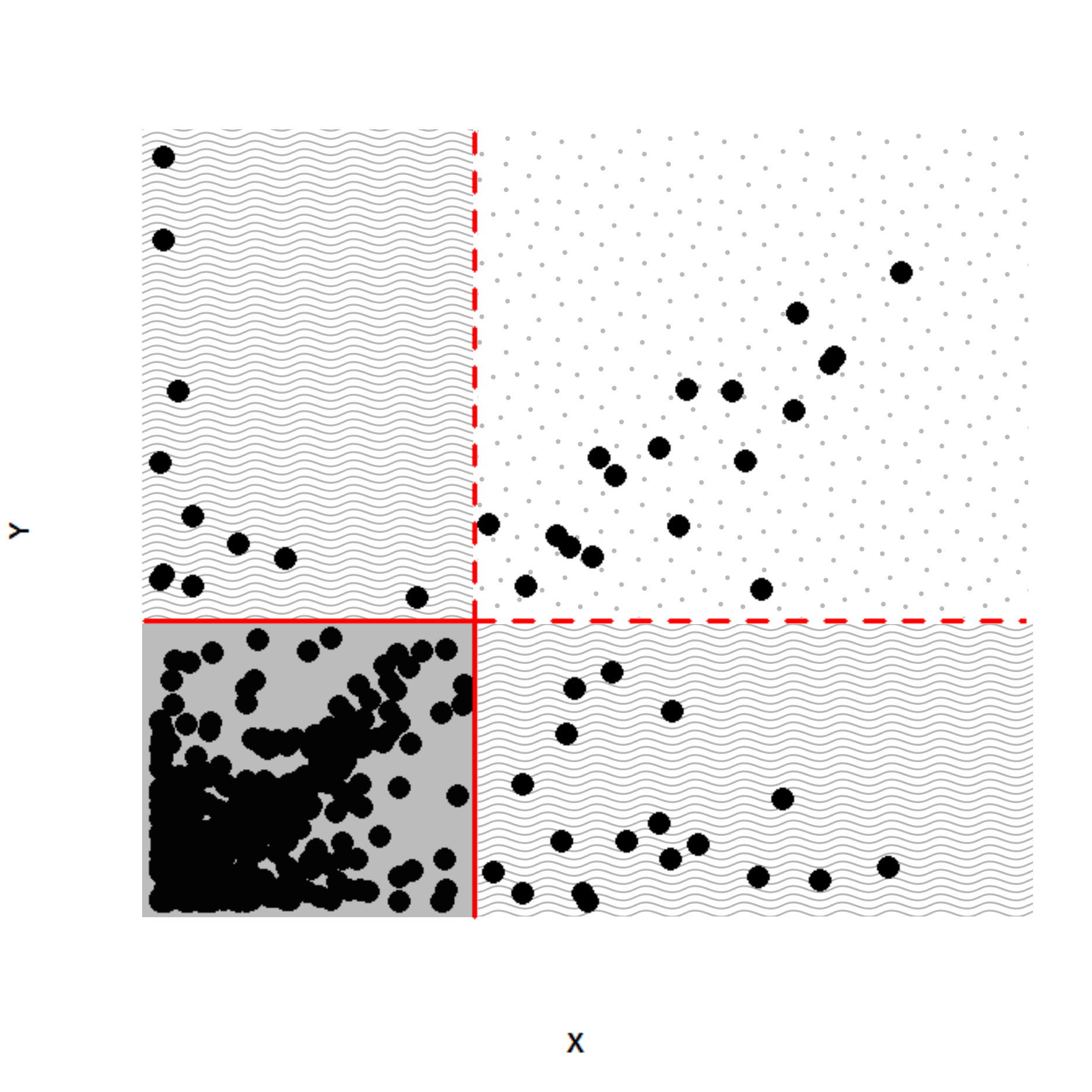}

\caption{Scatter plot of two asymptotically dependent variables.}
\label{importancesampling}
\end{figure}

Monte-Carlo simulations could be a convenient way of estimating the probability of joint occurrences of extreme events. However, as we recall from the introduction, the simulation of general multivariate heavy-tailed distributions is a serious issue. For instance, simulations of multivariate Generalised Pareto distributions can be carried out only in very specific cases as far as we know  (see \citet{Michel2007} for simulations in the logistic case). Nevertheless, in the particular case where we wish to estimate a probability of joint exceedances over a large threshold, the full simulation of the distribution over $\mathcal{O}$ is not needed. Mindful of the importance of sampling techniques, we propose to simulate the distribution over specific subspaces of  $\mathcal{O}$. The insight for the method is illustrated in Fig. \ref{importancesampling} where we have simulated two asymptotically dependent variables $X$ and $Y$. The figure emphasizes four regions but only the region with the dotted background is relevant if our interest lies in the probability of joint occurrences of large values.

We now introduce a novel procedure for estimating the probability of simultaneous exceedances over large thresholds. For illustrative purpose, we choose $p_m=(i_1,\dots,i_m)\in\mathcal{P}_d$ for some $m\leq d$ and fix a vector $\mathbf{x}=(x_1,\dots,x_m)$, each component being greater than $t$. We wish to estimate the probability of the set  $P_{m}(\mathbf{x}):=\left\{X^{(p_m)}>\mathbf{x}\right\}$.

By construction of the sets $E_i, \: i =1\dots l$, if there is no element $E\in\mathfrak{P}_d$ such that $p_m\subset E$,  then it will be inferred that the probability of $P_m(\mathbf{x})$ will be zero. Now, assume that there is one unique element $E\in\mathfrak{P}_d$ such that $p_m\subset E$ (generalisations when $E$ is not unique are straightforward). By construction, the probability of $P_m(t)$ is the same as the probability of $Q_m(t):=\left\{X^{(p_m)}>\mathbf{x},\:X^{(E\backslash p_m)}>t\right\}$. No assumption needs to be made regarding the components of $\bar{E}$.

Our estimation of the probability of $Q_m(t)$ is based on Eq. \eqref{limpolar} and uses the polar transformed data $(R,A):=T\left(X^{(E)}\right)$. Eq. \eqref{limpolar} states  that the angular and radial components $A$ and $R$ are asymptotically independent. Hence , assuming the angular and radial densities can be simulated, then the joint distribution can be simulated. For the simulation of the radial component, we assume its distribution is in the domain of attraction of an extreme value distribution (\cite{BGST}) so that a Generalized Pareto distribution can be fitted to its tail. In this paper the angular density was estimated with kernel estimators and was simulated with accept-reject methods.

Applying the inverse polar transform $T^{-1}$ to the simulated polar data, we obtain simulations of $X^{(E)}$ when all components are greater than $t$. The probability of $P_m(t)$ can be easily derived from these simulations. Note that the data are not necessarily identically distributed and are only assumed to have the same tail index. They can be rescaled to have the same order of magnitude by division by a high quantile of order $1-k/n$ for some $k$, $k\rightarrow\infty$, $k/n \rightarrow 0$.

\subsection{Numerical experiment}
\label{simu}

This paper deals primarily with two aspects of heavy-tail modelling and its application to risk assessment. Firstly, we showed in section \ref{asymptdepest} how the inference of the groups of asymptotically dependent variables made possible the estimation of the high dimensional angular probability, which had hitherto been intractable. Secondly, section \ref{jointExcProb} gave a new and efficient technique for estimating the probability of joint occurrence of extremal events. Therefore the simulation procedure needs to validate our clustering algorithm and then demonstrate the efficiency of the suggested heuristic method to estimate the probability of the joint occurrence of extreme events.\\

The simulation study in dimension $d=14$ is conducted as follows: we simulate $n$ realisations of a vector $X=\left(X^{(1)},\dots,X^{(14)}\right)$ of standard Pareto variables (so that $Z=X$). The dependence structure is modelled with a Gumbel copula  with dependence parameter $\nu \geq 1$ (\citet{Nelsen1999}), given by 

$$C_{\nu}\left(u_1,\dots,u_d\right)=\exp\left(-\left(\sum_{i=1}^{d}\left(-\log u_i\right)^{\nu}\right)^{1/\nu}\right).$$

The Gumbel copula accounts very efficiently for extremal dependences through its parameter $\nu$ (see \citet{Gudendorf2010}). The larger $\nu$, the more dependences there are, with asymptotic independence for $\nu=1$. We simulated five vectors with the following distributions

\begin{align*}
\left(X^{(1)},X^{(2)}\right) & \sim C_{\nu}\left(F\left(x^{(1)}\right),F\left(x^{(2)}\right)\right)\\
\left(X^{(3)},X^{(4)},X^{(5)}\right) & \sim \frac{1}{2}C_{\nu}\left(F\left(x^{(3)}\right),F\left(x^{(4)}\right)\right)F\left(x^{(5)}\right)+ \frac{1}{2}F\left(x^{(4)}),F(x^{(5)}\right)F\left(x^{(3)}\right) \\
\left(X^{(6)},X^{(7)},X^{(8)}\right) & \sim \frac{1}{2}C_{\nu}\left(F\left(x^{(6)}\right),F\left(x^{(7)}\right)\right)F\left(x^{(8)}\right)+ \frac{1}{2}F\left(x^{(7)}),F(x^{(8)}\right)F\left(x^{(6)}\right) \\
\left(X^{(9)},X^{(10)}\right) & \sim C_{\nu}\left(F\left(x^{(9)}\right),F\left(x^{(10)}\right)\right)\\
\left(X^{(11)},X^{(12)},X^{(13)},X^{(13)}\right) & \sim \frac{1}{2}C_{\nu}\left(F\left(x^{(11)}\right),F\left(x^{(13)}\right),F\left(x^{(14)}\right)\right)F\left(x^{(12)}\right)\\
                            &\qquad +{}\frac{1}{2}C_{\nu}\left(F\left(x^{(11)}\right),F\left(x^{(12)}\right)\right)F\left(x^{(13)}\right)F\left(x^{(14)}\right). \\
\end{align*}
where $\nu=2$ and $F(x) = 1-1/x$. This leads to the following 15 groups of asymptotically dependent variables:

\begin{itemize}
\item Singleton : $\{Z^{(3)}\}-\{Z^{(5)}\}-\{Z^{(6)}\}-\{Z^{(8)}\}-\{Z^{(12)}\}-\{Z^{(13)}\}-\{Z^{(14)}\}$.
\item Doublets : $\{Z^{(1)},Z^{(2)}\}-\{Z^{(3)},Z^{(4)}\}-\{Z^{(4)},Z^{(5)}\}-\{Z^{(6)},Z^{(7)}\}-\{Z^{(7)},Z^{(8)}\}\\-\{Z^{(9)},Z^{(10)}\}-\{Z^{(11)},Z^{(12)}\}$.
\item Triplets : $\{Z^{(11)},Z^{(13)},Z^{(14)}\}$.
\end{itemize}

Each group has the same weight $1/17$ except the groups $\{Z^{(1)},Z^{(2)}\}$ and $\{Z^{(9)},Z^{(10)}\}$ with weights $2/17$. We set $t=n/k$, where $k=k(n)\rightarrow\infty$ (see simulation results).

The complete results of our simulations are presented in table \ref{resSimu} where a type I error means that at least one good group was not discovered by the algorithm. A type II error means that at least one bad group was discovered. Most errors involved only one group meaning that only one good group was not discovered or only one bad group was discovered. The results show that the algorithm is very efficient even for small sample sizes.\\

\begin{table}
\label{resSimu}
\caption{Simulation results with  $e_j=0.2c_j$, $j=1\dots d$, $n_r=100$, $m_r=25$, $\sigma=0.05$.}
\fbox{%
\begin{tabular}{l|cccc}
 &  $\mathbf{n=1000}$  &  $\mathbf{n=2500}$ &   $\mathbf{n=5000}$ &  $\mathbf{n=10000}$ \\
  &  $\mathbf{k=100}$  &  $\mathbf{k=150}$ &   $\mathbf{k=250}$ &  $\mathbf{k=500}$ \\
 \hline
\textbf{No error}   & 76 & 90 & 94 & 99 \\
\textbf{Error I}    & 7  & 3  & 2  & 1\\
\textbf{Error II}   & 10 & 6  & 3  & 0\\
\textbf{Error I+II} & 7  & 1  &	1   &0
\end{tabular}}
\end{table}

\begin{figure}[!htb]
\centering
\includegraphics[width=100mm,height=90mm]{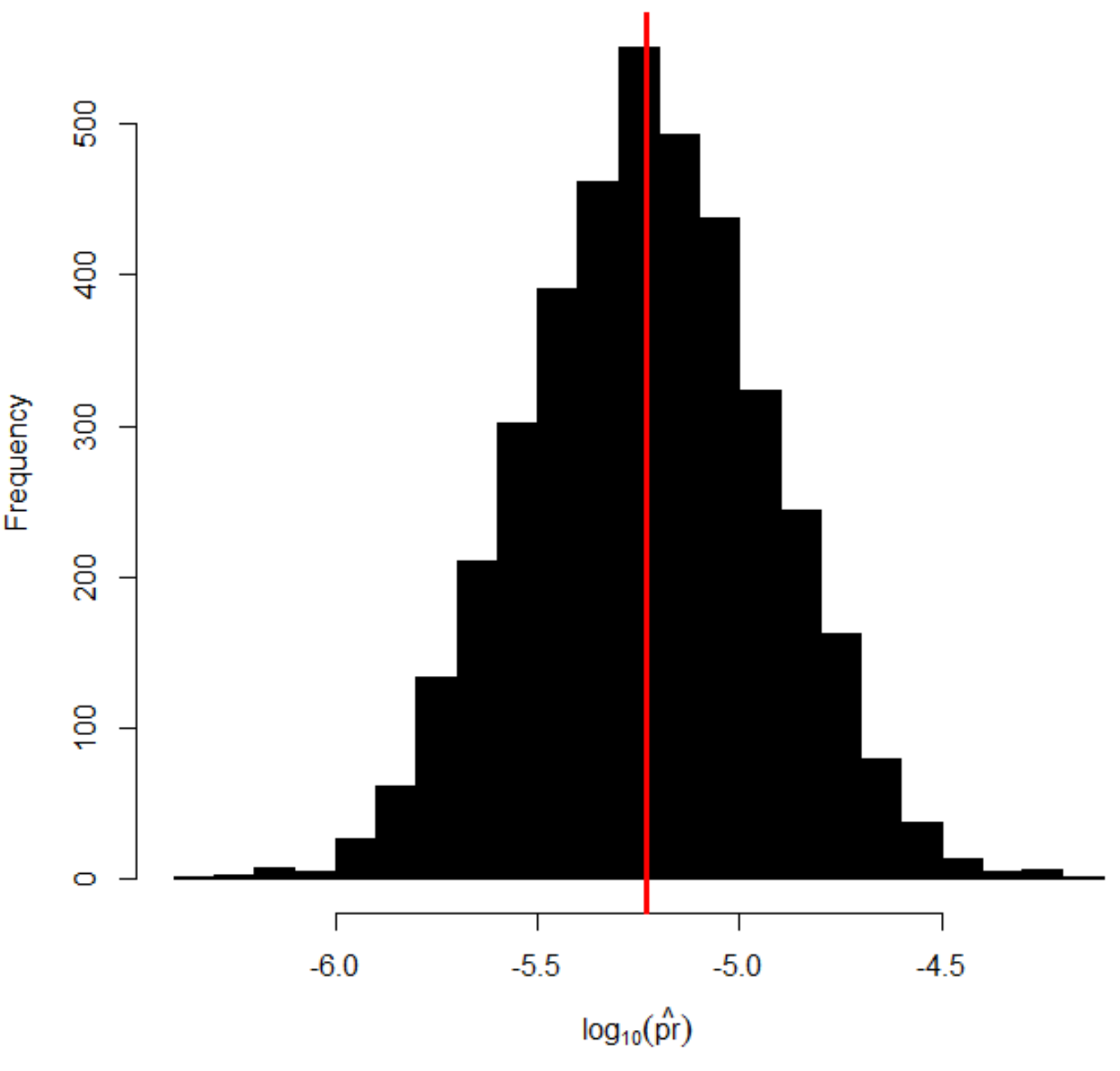}

\caption{Histogram of $\log_{10}\left(\widehat{P_2}(\mathbf{x})\right)$.  The value of $\log_{10}\left(P_2(\mathbf{x})\right)$ is indicated by the vertical solid line.}
\label{histoError}
\end{figure}

Now we set $n=10000$ and we wish to estimate the probability of $\Omega=\left\{ \left(X^{(1)},X^{(2)}\right) > \mathbf{x} \right\}$, where $\mathbf{x}=(100000,100000)$. We denote this probability by $P_2(\mathbf{x})$ and an estimator by $\widehat{P_2}(\mathbf{x})$. Note that the true value is $P_2(\mathbf{x})=5.86\times 10^{-6}$  meaning that $\Omega$ is only observed once out of 17 samples of size $n$. we repeat the experiment $N=10000$ times and plot the histogram of $\log_{10}\left(\widehat{P_2}(\mathbf{x})\right)$ in Fig. \ref{histoError}. The mean relative error is 0.046 which means that $P_2(\mathbf{x})$ is over or under estimated by a factor of $11.1\%$ on average.

\section{Case study: Risk assessment in the sloshing industry}
\label{caseStudy}
We now apply the methodology introduced in the first part of this paper to assess the risk associated with the sloshing phenomenon in the LNG shipping industry.

\subsection{Assessing groups of asymptotically dependent sensors}
\label{clustering}
\begin{figure}[!htb]
\centering
\includegraphics[width=100mm,height=90mm]{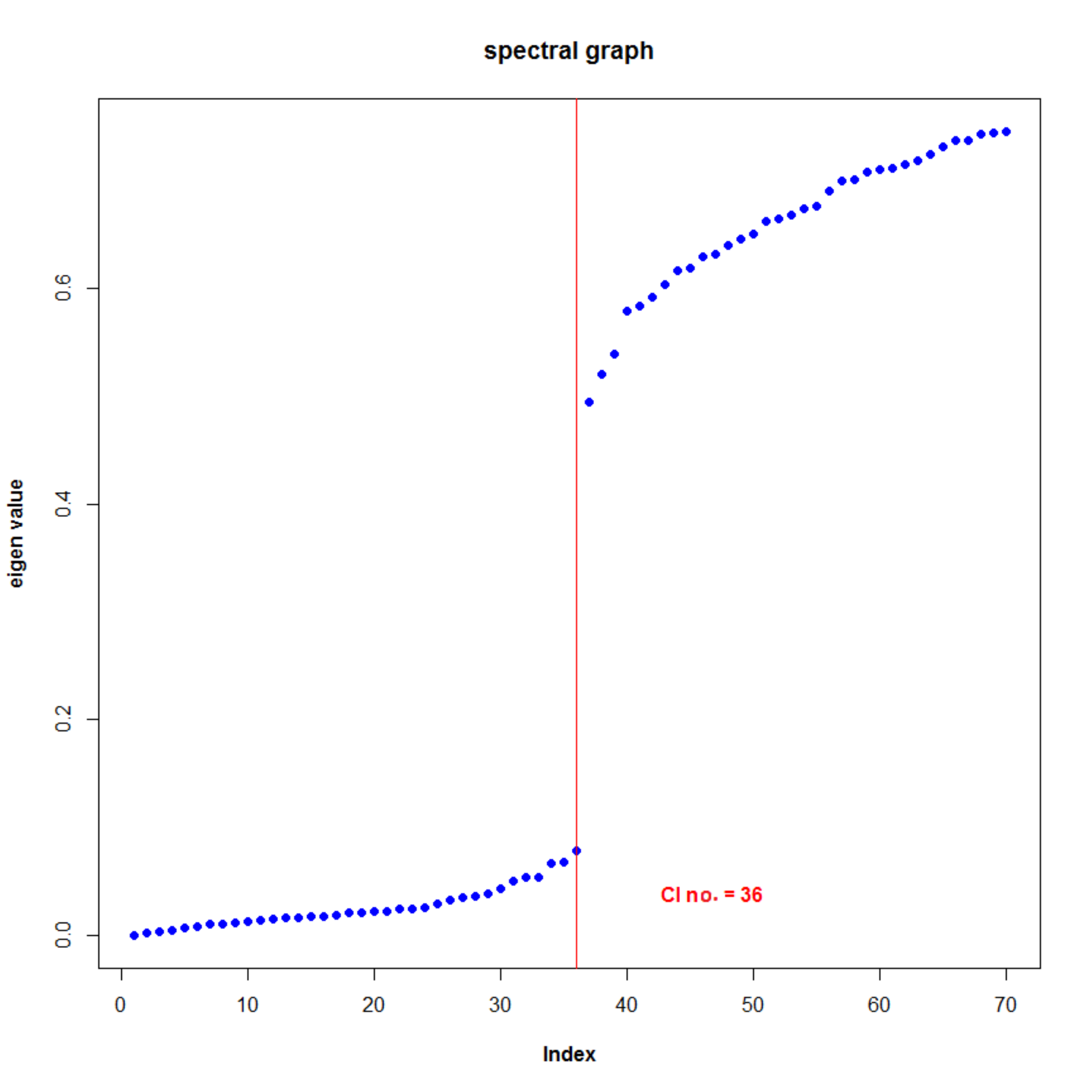}

\caption{Spectral graph for the sloshing data set.}
\label{SpectralGraph}
\end{figure}

\begin{figure}[!htb]
\centering
\includegraphics[width=147mm,height=96mm]{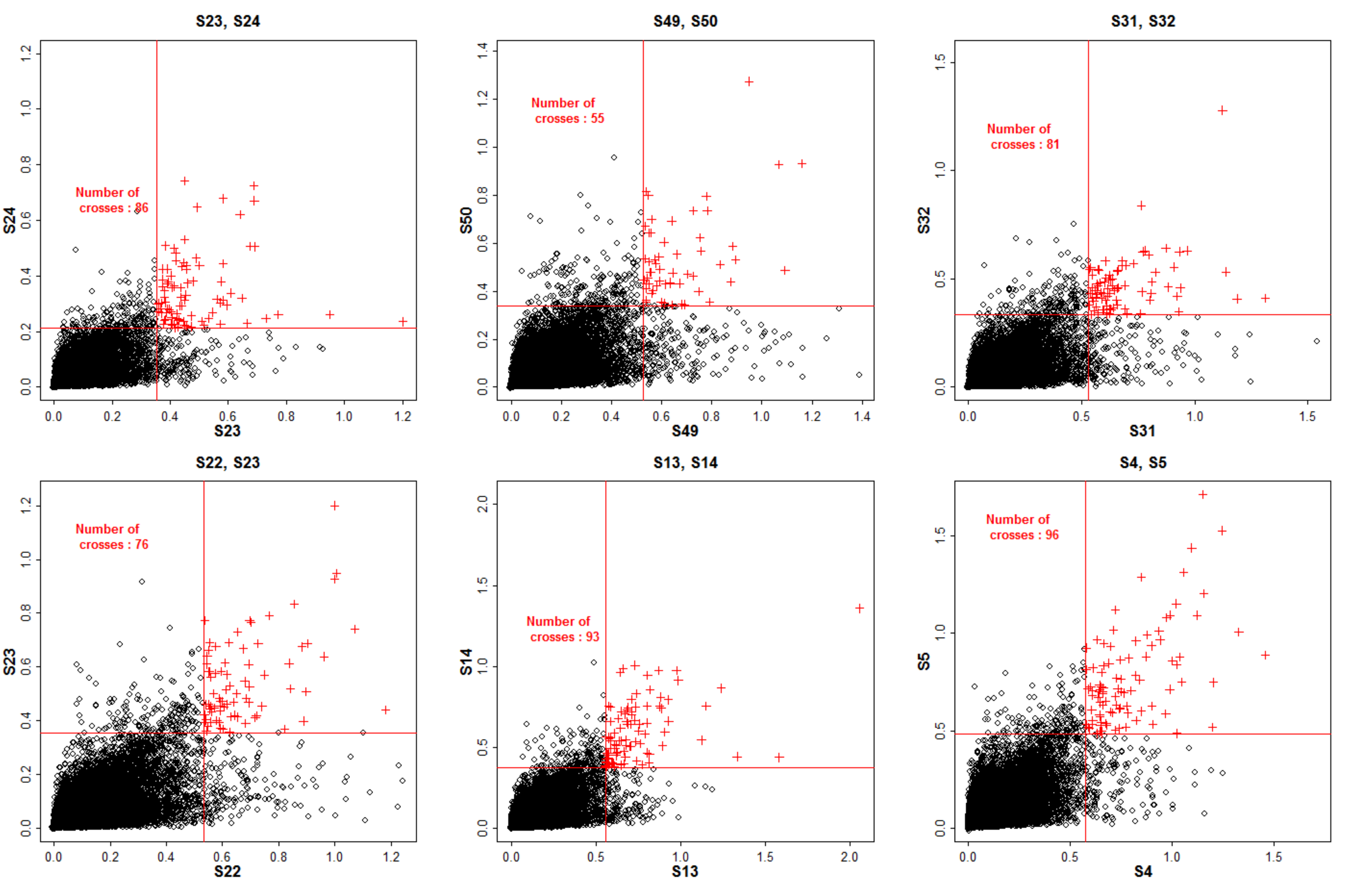}

\caption{Scatter plots for the two-dimensional clusters. The corresponding quantile of order $1-k/n$ is emphasized by the horizontal and vertical solid lines. The crosses represent the impacts where both sensors exceeded their quantile.}
\label{scatterExtr}
\end{figure}

The extremal spectral clustering algorithm is used to estimate the groups of asymptotically dependent sensors in the sloshing data set with the following parameters: $k = 250$, $n_r=100$, $m_r=50$, $e_i=0.25c_i$ for any cluster $C_i$ of size $c_i$. As with the simulation study, we consider the standardised version of the data set. Fig. \ref{SpectralGraph}, in which the first 70 eigenvalues are plotted, strongly advocates for the existence of 36 clusters. The results are displayed in Table \ref{SEC}: the data exhibit few asymptotic dependences, most clusters being singletons and the largest groups having dimension 2. This was somewhat predictable, insofar as most phenomena characterising sloshing are very local, being typically the size of one sensor. Notice from Table \ref{HillCI} that we cannot reject the hypothesis that two sensors belonging to the same group of dimension 2 have the same tail-index and then; in what follows we then assume that we are in the so-called  \textit{standard case}. In Fig. \ref{scatterExtr} we draw the scatter plot of the pressure measurements of any of the 2-dimensional groups. It shows that these sensors clearly exhibit asymptotic dependences. The results of the estimation of $\pi_p$ are also presented in Table \ref{SEC}.

\begin{table}
\caption{\label{SEC}Results of the Extremal Spectral Clustering}
\begin{tabular}{ccc}
\textbf{Clusters} & \textbf{Frequencies} (\%) & $\pi_p$ \\
 \hline 
$\{S8\}$---$\{S54\}$---$\{S53\}$& \multirow{9}{*}{ranging from 73 to 98 } &$3.7\times 10^{-2}$---$4.2\times 10^{-2}$---$3.2\times 10^{-2} $\\
$\{S45\}$---$\{S40\}$---$\{S18\}$ & & $3.4\times 10^{-2}$---$4.1\times 10^{-2}$---$3.5\times 10^{-2} $\\
$\{S26\}$---$\{S35\}$---$\{S43\}$ & & $2.6\times 10^{-2}$---$2.4\times 10^{-2}$---$2.7\times 10^{-2} $\\
$\{S36\}$---$\{S16\}$---$\{S7\}$ & & $3.0\times 10^{-2}$---$2.1\times 10^{-2}$---$3.2\times 10^{-2} $\\
$\{S52\}$---$\{S9\}$---$\{S44\}$ & & $3.0\times 10^{-2}$---$4.7\times 10^{-2}$---$3.4\times 10^{-2} $\\
$\{S27\}$---$\{S4\}$---$\{S49\}$ & & $3.1\times 10^{-2}$---$3.2\times 10^{-2}$---$4.8\times 10^{-2} $\\
$\{S22\}$---$\{S51\}$---$\{S41\}$ & & $3.4\times 10^{-2}$---$3.2\times 10^{-2}$---$2.7\times 10^{-2} $\\
$\{S34\}$---$\{S6\}$---$\{S31\}$ & & $2.0\times 10^{-2}$---$3.4\times 10^{-2}$---$3.6\times 10^{-2} $\\
$\{S42\}$---$\{S25\}$---$\{S13\}$ & & $2.1\times 10^{-2}$---$1.9\times 10^{-2}$---$3.2\times 10^{-2} $\\
$\{S17\}$ & & $2.8\times 10^{-2}$\\
\hline 
\hline 
$\{S23, S24\}$ & 86 & $1.8\times 10^{-2} $\\
$\{S49, S50\}$ & 92 & $1.0\times 10^{-2} $\\
$\{S31, S32\}$ & 89 & $1.9\times 10^{-2} $\\
$\{S22, S23\}$ & 83 & $1.7\times 10^{-2} $\\
$\{S13, S14\}$ & 92 & $2.4\times 10^{-2} $\\
$\{S4, S5\}$ & 93 & $2.5\times 10^{-2} $\\
\hline 

\end{tabular}
\end{table}

\subsection{Application to the estimation of the joint occurrences of high pressures on several sensors}
GTT designs its vessels so that the probability of failure of the cargo containment system is less than a target probability of $10^{-3}$ in forty years (recall that the small scale data set corresponds to 6 months at full scale). A failure occurs if the pressure loads exerted on the membrane are too large, and hence the  areas most likely to be exposed to such loads need to be reinforced. The maximal admissible load is a function of the impacted area. According to the dependence structure identified in section \ref{clustering} at most two sensors can be impacted at the same time by large pressure loads. For an area the size of two sensors ((1 cm${}^2$)), this pressure is approximately 1.5 bar. Note that in one tank, there are 12 arrays with identical behaviour and there are four tanks in the vessel. Hence the probability that in one array the pressure exerted on an area of 1 cm${}^2$ is greater than 1.5 bar needs to be multiplied by 48 to obtain the equivalent probability for one tank.The purpose of the remainder of this section is to estimate this bivariate probability for the sensor array.

\begin{figure}[!htb]
\centering
\includegraphics[width=147mm,height=125mm]{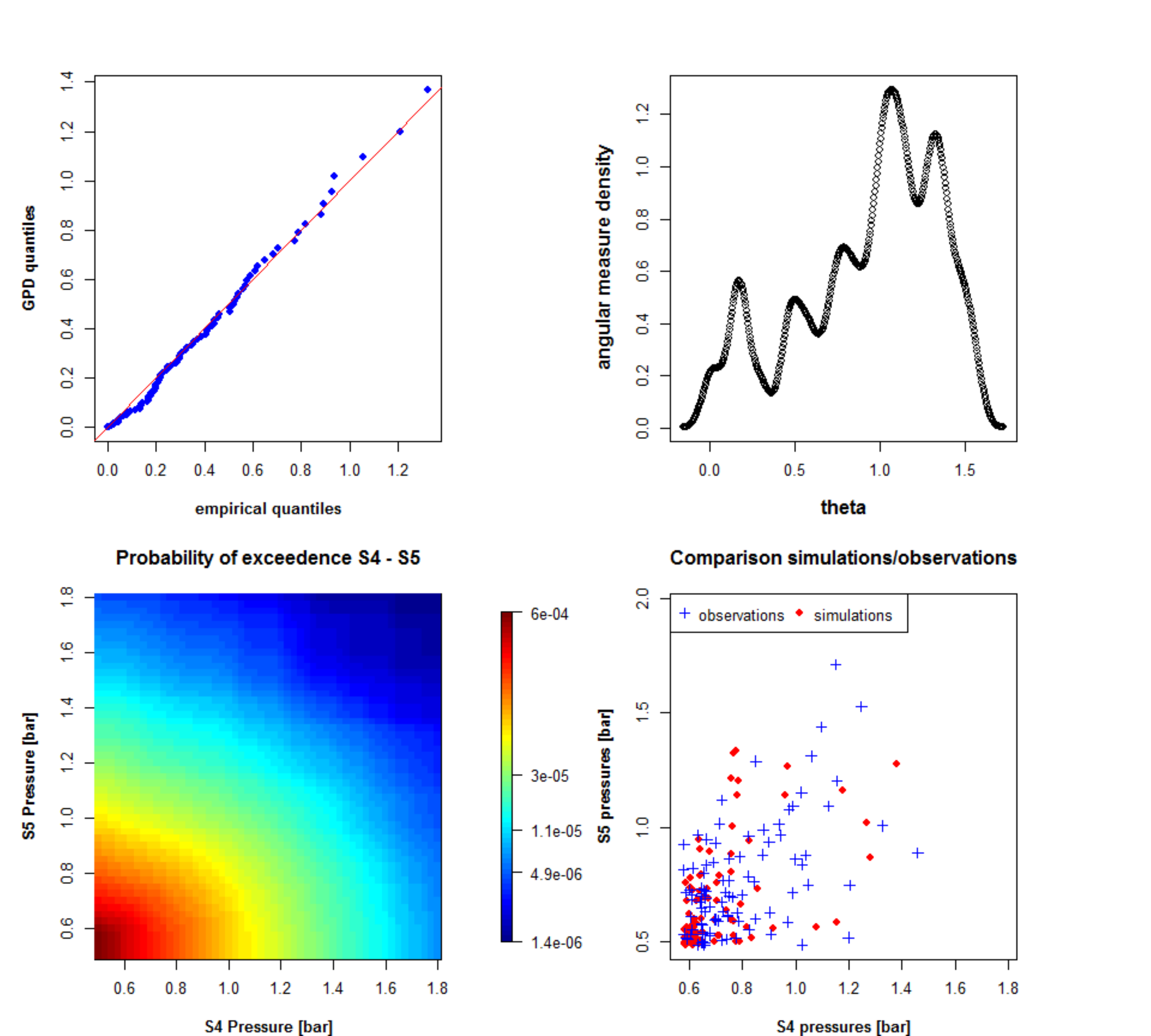}

\caption{Estimation of the joint probability of exceedance for sensors $S4$ and $S5$. Top-left: GPD fit to the radius. Top-right: estimation of the angular measure. Bottom-left: Estimation of the probability of bivariate exceedance. Bottom-right: comparison between the observed data and the simulated data, with same sample size.}
\label{trueDataEst}
\end{figure}

\begin{figure}[!htb]
\centering
\includegraphics[width=75mm,height=75mm]{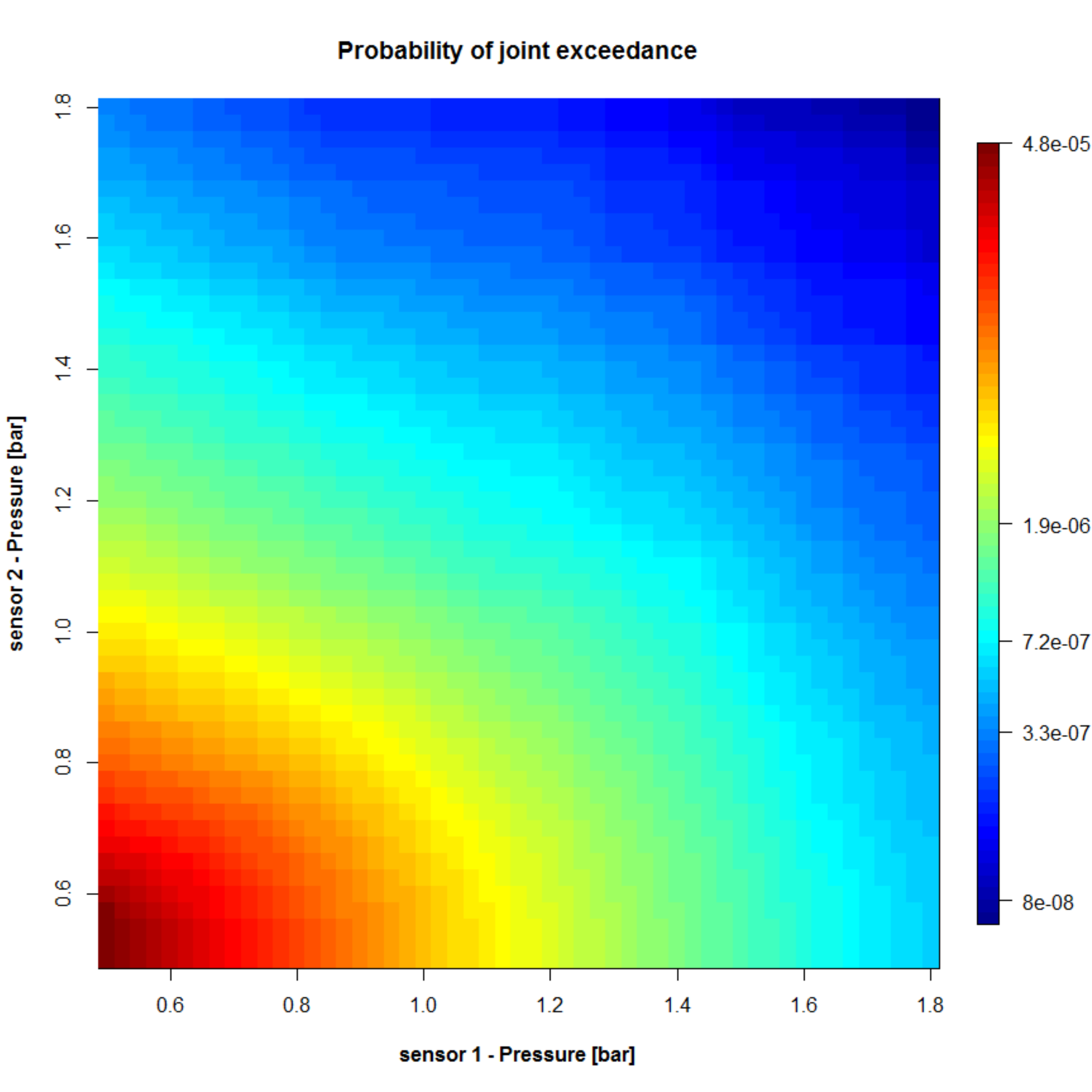}

\caption{Estimation of the probability that two sensors exceed a large pressure value.}
\label{globalexcessproba}
\end{figure}

The complete procedure (that is, the GPD fit to the radial component and estimation of the angular density) for the estimation of the joint exceedance is detailed in figure \ref{trueDataEst} for sensors $S4$ and $S5$. The overall result for the array, i.e. the probability that two sensors jointly record large values, is given in Fig. \ref{globalexcessproba}. Results group by group are also provided in Fig.\ref{jointexcess}, appendix \ref{statdesc}. In Table \ref{exceedenceProb}, we focus on the particular case of exceedance greater than 1.5 bar. The overall probability that two sensors simultaneously exceed 1.5 bar over forty years in the tank is $1.28\times 10^{-5}$.

\begin{table}
\caption{\label{exceedenceProb}Probability of bivariate exceedance by group. The value $\hat{p}$ stand for the estimation of the probability of simultaneous exceedance over 1.5 bars for the two sensors of the group.}
\fbox{%
\begin{tabular}{lllllll}
 & $\{S23-S24\}$ & $\{S49-S50\}$ & $\{S31-S32\}$ & $\{S22-S23\}$ & $\{S13-S14\}$ & $\{S4-S5\}$ \\
$\hat{p}$ & 0 & $4.3\times10^{-7}$ & $1.7\times10^{-7}$ & $3.4\times10^{-7}$ & $1.8\times10^{-6}$ & $3.8\times10^{-6}$\\
$\pi_i$ & $3.1\times10^{-2}$ & $2.4\times10^{-2}$ & $2.9\times10^{-2}$ & $2.9\times10^{-2}$ & $3.2\times10^{-2}$ & $3.6\times10^{-2}$ \\
 $\pi_i\hat{p}$ & 0 & $1.0\times10^{-8}$ & $5.0\times10^{-9}$ & $9.9\times10^{-9}$ & $5.7\times10^{-8}$ & $1.4\times10^{-7}$ \\
 
\end{tabular}}
\end{table}

\section{Discussion and outlook for the future}
\label{discussion}
A very high pressure is fortunately an extreme, and rare, event and it seemed reasonable to investigate the joint distribution of such pressures through heavy-tail analysis. This is a very common and conservative approach in risk assessment because it is unlikely that it leads to an underestimate of the risk. Our goal was to study non-parametrically the extremal dependencies among observed extremal pressures and to estimate the probability of simultaneous occurrences of large pressures at different locations in the tank. This was not possible directly because the dimension of the data set we consider is very high (in the example considered, the dimension is 36). So far, classical methods cannot deal with more than three or four dimensions. To overcome this issue, we proposed a novel latent variable analysis of the angular measure that enabled us to overcome the 'curse of dimensionality' and render its estimation tractable even in large dimensions. This major breakthrough makes multivariate heavy-tail modelling possible, even for high dimensional data sets. 

The statistical techniques proposed in this paper showed their capacity to exhibit groups of asymptotically dependent sensors in the simulation experiments we carried out. Our approach makes hitherto intractable multivariate risk analysis possible. We provide a method for estimating the probability of the simultaneous exceedance of a high threshold of the pressures recorded by the sensors.

Several tuning parameters may have a large influence on the results and the cooperation with GTT's sloshing expert was of great value. The first parameter is the number $k$ of extremes used. Its choice is always a trade-off between bias and variance and, in this paper, a result also of physical considerations with experts wishing to focus on the largest pressure peaks. Second, the parameter $\sigma$ of the similarity function used to compute the graph Laplacian can have a dramatic influence, though the optimal number of 36 clusters seemed quiet clear. A sensitivity study was conducted and the results did not change for wide ranges of $\sigma$. In the end, the most influential parameter seemed to be the threshold $e_i$, designed to control the size of the clusters. A choice of smaller thresholds $e_i$ may have led to the discovery of larger clusters. However, we point out that it is common in real data sets for variables to exhibit few asymptotic dependences and therefore asymptotic independences are frequent. Furthermore, it is known by sloshing experts that sloshing pressure peaks are sharp and it seemed reasonable for our applications to consider sufficiently large thresholds $e_i$ to avoid the inclusion of asymptotically independent sensors in our groups of asymptotically dependent variables. We emphasize the fact that the proposed methodology is very general and can be used for any multidimensional heavy-tailed data set quite apart from the specific case of sloshing data.

\section{Acknowledgement} The authors are very grateful to GTT engineers for their help and invaluable advice.

\bibliographystyle{chicago}
\bibliography{bibliographieAbrev}

\newpage

\begin{appendix}

\section{Descriptive statistics}
\label{statdesc}
\begin{table}
\caption{\label{quantilesallsens}High quantiles for all the sensors of the array}
\fbox{%
\begin{tabular}{c|ccccc}
 & \multicolumn{4}{c}{\textbf{quantiles}}  \\
\textbf{Sensor}  & $\mathbf{0.9}$  & $\mathbf{0.99}$  & $\mathbf{0.999}$  & $\mathbf{0.9999}$  & \textbf{max} \\
\hline \\
\textbf{S13}  & 0.099  & 0.312  & 0.623  & 0.929  & 2.061 \\
\textbf{S5}  & 0.101  & 0.268  & 0.554  & 0.977  & 1.712 \\
\textbf{S31}  & 0.092  & 0.293  & 0.603  & 0.987  & 1.542 \\
\textbf{S4}  & 0.112  & 0.327  & 0.664  & 1.067  & 1.459 \\
\textbf{S49}  & 0.092  & 0.294  & 0.609  & 1.020  & 1.391 \\
\textbf{S40}  & 0.090  & 0.286  & 0.582  & 0.934  & 1.360 \\
\textbf{S14}  & 0.074  & 0.203  & 0.444  & 0.794  & 1.358 \\
\textbf{S41}  & 0.064  & 0.185  & 0.388  & 0.680  & 1.321 \\
\textbf{S32}  & 0.063  & 0.181  & 0.393  & 0.588  & 1.279 \\
\textbf{S50}  & 0.065  & 0.190  & 0.409  & 0.695  & 1.274 \\
\textbf{S22}  & 0.092  & 0.301  & 0.605  & 1.000  & 1.242 \\
\textbf{S23}  & 0.067  & 0.188  & 0.418  & 0.687  & 1.201 \\
\textbf{S6}  & 0.054  & 0.132  & 0.277  & 0.477  & 1.033 \\
\textbf{S15}  & 0.056  & 0.128  & 0.273  & 0.532  & 0.837 \\
\textbf{S24}  & 0.051  & 0.119  & 0.246  & 0.458  & 0.743 \\
\textbf{S33}  & 0.050  & 0.117  & 0.253  & 0.426  & 0.743 \\
\textbf{S7}  & 0.056  & 0.128  & 0.254  & 0.418  & 0.715 \\
\textbf{S16}  & 0.050  & 0.101  & 0.181  & 0.326  & 0.697 \\
\textbf{S8}  & 0.048  & 0.109  & 0.206  & 0.327  & 0.673 \\
\textbf{S51}  & 0.044  & 0.106  & 0.214  & 0.364  & 0.665 \\
\textbf{S42}  & 0.049  & 0.115  & 0.233  & 0.390  & 0.660 \\
\textbf{S25}  & 0.044  & 0.091  & 0.166  & 0.301  & 0.642 \\
\textbf{S43}  & 0.040  & 0.087  & 0.158  & 0.263  & 0.572 \\
\textbf{S9}  & 0.051  & 0.113  & 0.217  & 0.369  & 0.561 \\
\textbf{S18}  & 0.039  & 0.081  & 0.139  & 0.225  & 0.560 \\
\textbf{S34}  & 0.042  & 0.089  & 0.167  & 0.279  & 0.560 \\
\textbf{S35}  & 0.037  & 0.076  & 0.138  & 0.227  & 0.434 \\
\textbf{S26}  & 0.041  & 0.082  & 0.142  & 0.230  & 0.390 \\
\textbf{S52}  & 0.036  & 0.079  & 0.147  & 0.248  & 0.390 \\
\textbf{S17}  & 0.051  & 0.093  & 0.156  & 0.263  & 0.388 \\
\textbf{S53}  & 0.032  & 0.072  & 0.124  & 0.187  & 0.332 \\
\textbf{S36}  & 0.033  & 0.071  & 0.124  & 0.196  & 0.319 \\
\textbf{S27}  & 0.037  & 0.075  & 0.129  & 0.201  & 0.301 \\
\textbf{S44}  & 0.034  & 0.073  & 0.128  & 0.214  & 0.299 \\
\textbf{S45}  & 0.032  & 0.070  & 0.126  & 0.201  & 0.296 \\
\textbf{S54}  & 0.030  & 0.069  & 0.124  & 0.194  & 0.226 \\

\end{tabular}}
\end{table}

\begin{figure}[!htb]
\centering
\includegraphics[width=90mm,height=70mm]{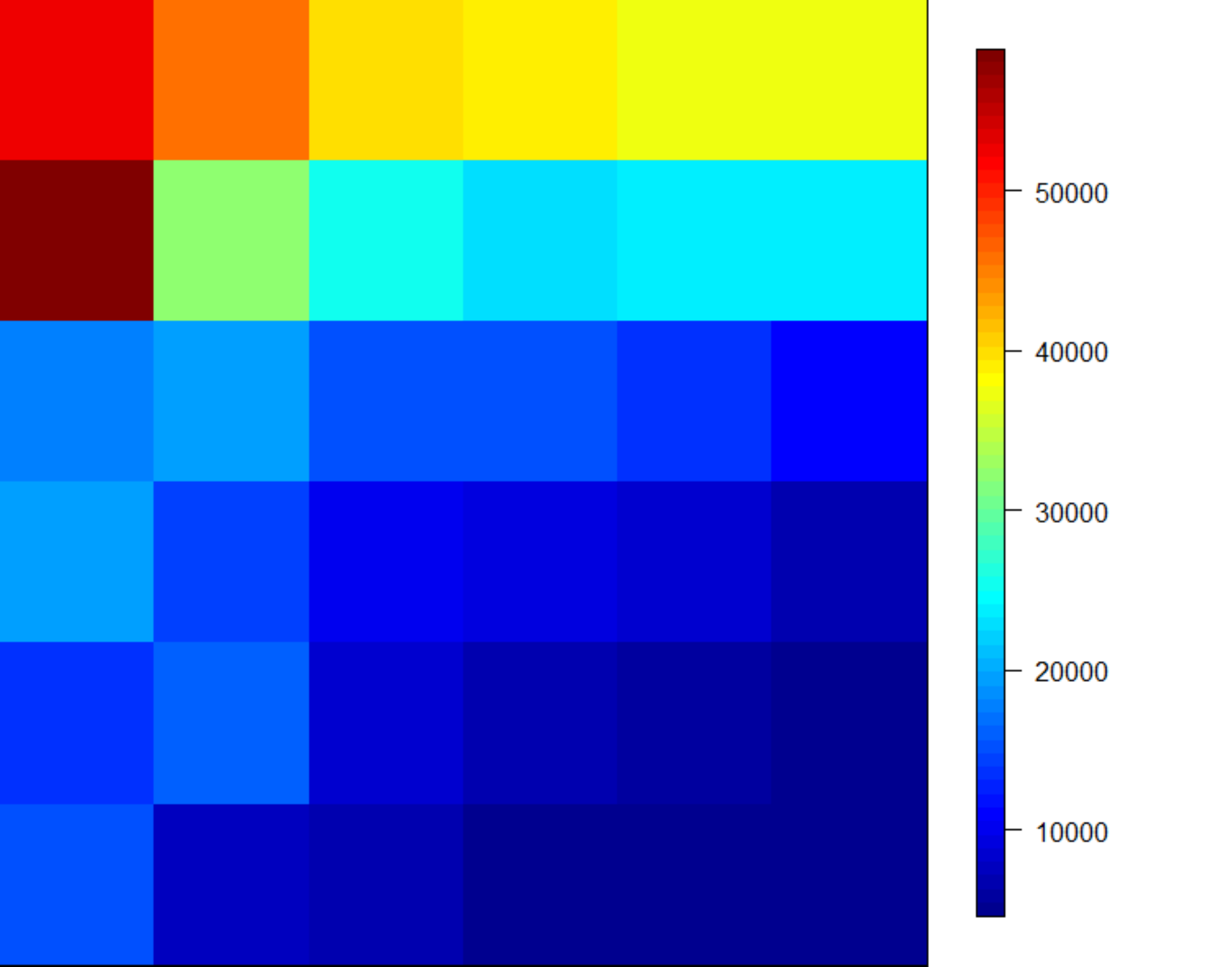}

\caption{Spatial map of the number of impacts detected along the array.}
\label{nb_evt_sens}
\end{figure}

\section{Spectral Clustering algorithm}
\label{SpectralClustering}
\begin{figure}[!htb]
\centering
\includegraphics[width=147mm,height=90mm]{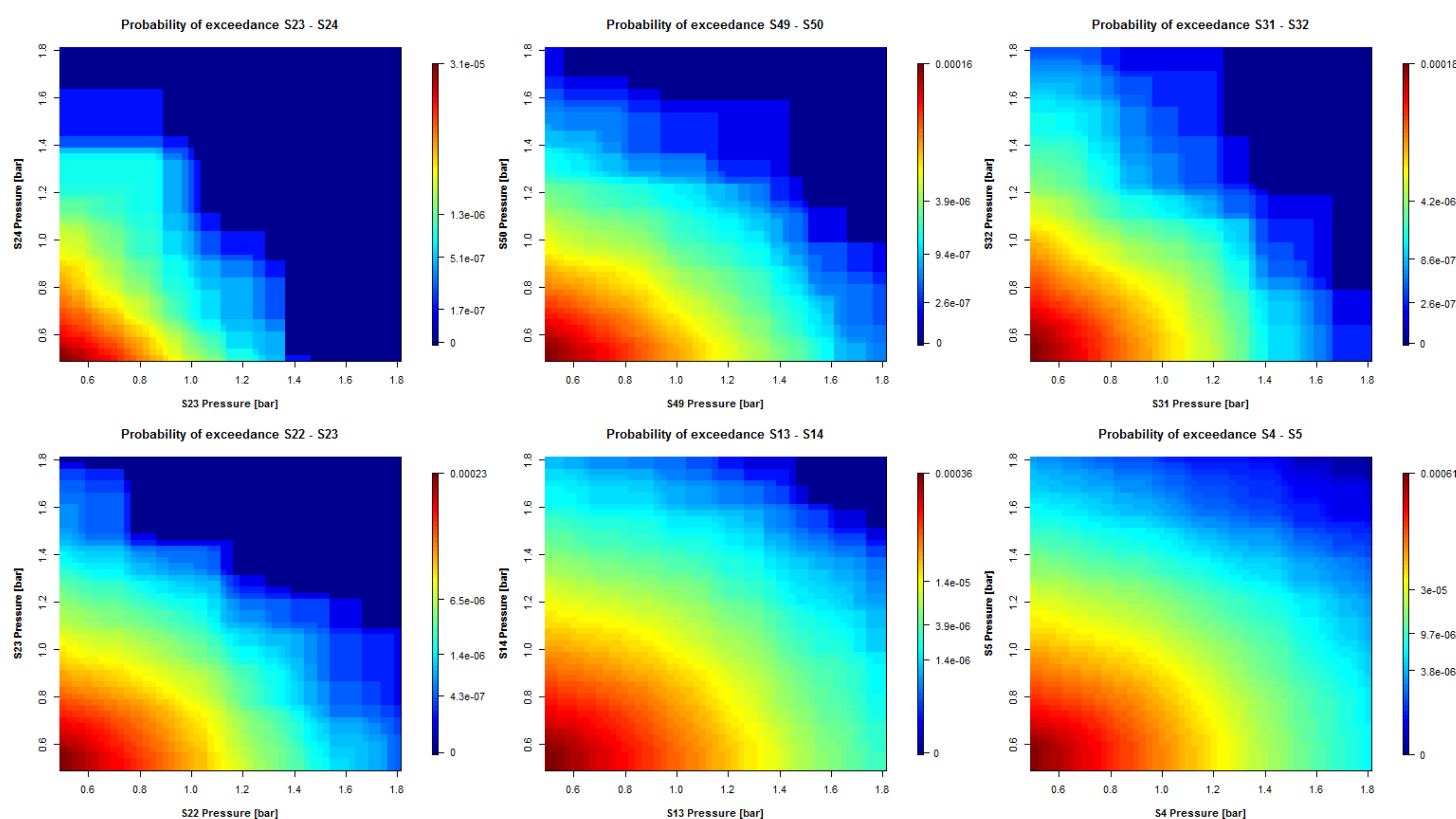}

\caption{Estimation of the probability of simultaneous occurrence of large pressures for each group. }
\label{jointexcess}
\end{figure}

\begin{center}
\fbox{
\begin{minipage}{0.8\textwidth}
\medskip

\centering
\textbf{\textit{Normalised Spectral Clustering}}\\
\begin{itemize}
\item[]\textbf{Input}: Similarity matrix $W$.
\item Build  the similarity graph $(V,E)$ with weighted adjacency matrix $W$.
\item Compute the normalised Laplacian $L_{sym}$ and let $k$ be the dimensionality of the eigenvalue $0$.
\item Compute $k$ orthonormal eigenvectors $t_1,\; \ldots,\; t_k$ of $L_{sym}$ and let $T \in \R^{n\times k}$ be the matrix with vectors $t_1,\; \ldots,\; t_k$ as columns.
\item For $i = 1,\dots , n$, let $y_i \in \R^{k}$ be the vector corresponding to the $i^{th}$ row of $T$. Segment the set of points $\{y_i:\;\; i=1,\;\ldots,\; n\}$ into clusters $C_1,\; \ldots,\; C_k$ using the $k$-means algorithm.
\item[]\textbf{Output}: $C_1,\; \ldots,\; C_k$.
\begin{flushleft}
\end{flushleft}
\end{itemize}

\end{minipage}
}
\end{center}
\end{appendix}

\end{document}